\documentclass[final]{siamltex}
\usepackage{epsfig,float,amssymb,amsmath,tabularx,lineno}
\newcommand{\pt}{\,\partial_t\,}
\newcommand{\px}{\,\partial_x\,}
\newcommand{\pz}{\,\partial_z\,}
\newcommand{\vq}{\vec q}
\newcommand{\A}{\,\mathcal A\,}
\newcommand{\B}{\mathcal B\,}
\newcommand{\C}{\,\mathcal C\,}
\renewcommand{\S}{\mathcal S\,}

\newcommand{\dt}{\Delta t }
\newcommand{\dx}{\Delta x }
\newcommand{\dz}{h }
\newcommand{\lb}{\langle\, }
\newcommand{\rb}{\,\rangle }

\title{High-order finite volume schemes for layered atmospheric models}

\author{Dante Kalise\thanks{Department of Mathematics, University of Bergen, Norway ({\tt dante.kalise@math.uib.no}).}
        \and Ivar Lie\thanks{ R \& D Department, StormGeo AS, Oslo, Norway ({\tt ivli@online.no}).}
        \and Eleuterio F. Toro \thanks{Laboratory of Applied Mathematics, Department of Civil and Environmental Engineering, University of Trento, Italy ({\tt toro@ing.unitn.it}).}}

\begin{document}

\maketitle

\begin{abstract}
We present a numerical scheme for the solution of a class of atmospheric models where high horizontal resolution is required while a coarser vertical structure is allowed. The proposed scheme considers a layering procedure for the original set of equations, and the use of high-order ADER finite volume schemes for the solution of the system of balance laws arising from the dimensional reduction procedure. We present several types of layering based upon Galerkin discretizations of the vertical structure, and we study the effect of incrementing the order of horizontal approximation. Numerical experiments for the computational validation of the convergence of the scheme together with the study of physical phenomena are performed over 2D linear advective models, including a set of equations for an isothermal atmosphere.
\end{abstract}



\pagestyle{myheadings}
\thispagestyle{plain}
\markboth{KALISE, LIE AND TORO}{HIGH-ORDER FV SCHEMES FOR LAYERED ATMOSPHERIC MODELS}

\section{Introduction}

When we study phenomena in the atmosphere, they are often confined to
certain layers, for example the atmospheric boundary layer (ABL) and the
stratosphere. Models of atmospheric phenomena often reflects this layering.
Examples of use of layered models are: air pollution models (see for example
the early study in \cite{reible}), moisture and precipitation studies
\cite{lambaerts}, hurricane modelling \cite{lam-10},
large scale atmospheric dynamics and climate models.

There are  differences between a conventional vertical discretization of a 3D
atmospheric model, and a layered model. The most important difference is that
layered models are focused on specific properties of certain layers,
whereas the purpose of a conventional vertical discretization is to express
numerically the vertical variation of the complete atmosphere.
Vertical discretization do treat say the ABL specifically, but only by
having sufficiently man vertical layers so as to resolve the phenomena
of interest.

The numerical schemes often used in atmospheric models are some form of
conservative difference schemes in the vertical and second order centered
difference schemes on a staggered grid in the horizontal.
The global weather model at the European Centre for Medium Range Forecasts (ECMWF)
introduced a finite element scheme in the vertical relatively recently,
and finite volume schemes has also been constructed \cite{ahil-vert}.

What we want to communicate here is an investigation of using of a high-
order ($\geq 2$) finite volume discretization  in the horizontal combined with
several different schemes in the vertical. The purpose is to find out
how efficient and accurate such discretizations are.
The horizontal discretization that we apply is a finite-volume high-order Godunov
scheme called ADER (Arbitrary high-order DERivatives), proposed by Toro and Titarev (see \cite{ader-1} and \cite{ader-2}, among others). The ADER scheme is actually a family of schemes with increasing order. Hence there is a trade-off
between order and computational complexity.
The schemes have been shown to work well in a series of applications,
see for example \cite{dumbser} and \cite{ader-3d}; related schemes, applied to shallow water models, are presented in \cite{lam-28}, \cite{lam-30} and \cite{lam-31}.

A layered model is obtained by some form of averaging of the governing
equations in the specified layer. This procedure then leads to
a set of equations in the horizontal, for each layer; the equations in the layers are then coupled. An example of such procedure for the linearized 3D Euler equations with an orthogonal set of finite element hat functions can be seen in \cite{temam}, obtaining what is called a ``2.5D'' or ``1.5D'' model.
We extend this approach by not only considering piecewise linear finite elements, but also a pseudosprectal method, a layering based on discontinuous Galerkin methods and a method based on vertical averaging imposing conservation properties. Such investigation in the layering procedure arises from the fact that apart from the physical basis of a layering, there is a potential saving in computational time, since only a few layers are required
and application of numerical schemes specifically aimed at 2D or 1D
problems can be used.

One of the motivations for our work is its application to data
assimilation, more specifically ensemble-type of methods as the Ensemble Kalman filter (EnKF).
For an overview, see e.g. \cite{hamill}.
The EnKF procedure involves an ensemble of model runs between each
observation time point, and preferable the size of this ensemble should
be as large as practically possible. This implies that model runs should
be as efficient as possible, and thus, the use of dimensionally reduced models is highly desirable.
It is worth to emphasize that our focus on data assimilation is somewhat
different from conventional data assimilation in weather models
( treated in e.g. \cite{kalnay} and \cite{daley}).
We are concentrating on high resolution
models where there are a few observations in the ABL, and we want to
study the effect of these observations in the other layers of the
atmosphere. A typical application where such a situation is relevant,
is windpower modelling. See \cite{ahil-enkf} for a detailed description.

To computationally study convergence, and to assess the performance of the method,
we first use a linear 2D advection where we investigate
the properties of the numerical schemes. Then we proceed to approximate a linear, non-hydrostatic
atmospheric model with the same family of schemes; apart from convergence issues, we present a physical test related to the propagation of inertia-gravity waves of non-hydrostatic scale in a periodic channel.

This paper is organized as follows: In section 2 we describe in detail the
verical discretization (layering) methods for a  2D linear advection-reaction problem. Section 3 contains a brief description of the ADER scheme for horizontal discretization. Numerical tests validating convergence rates and other relevant physical phenomena for 2D linear advection and for a 2D linear atmospheric model are shown. Outlook and summary is included in section 5.

\section{Vertical discretization}
We consider a linear, constant-coefficient, two-dimensional advection-reaction model given by

\begin{equation}\label{ley}
\pt q +a \px q+b\pz q+s \,q=0.
\end{equation}
The problem is spatially defined over a rectangular domain $\Omega=\{(x,z)\in[0,L_x]\times[0,H_z]\}$ together with compatible boundary and initial conditions.We study four different methods to reduce the dimension of the model by discretizing the vertical component. It is relevant to note at this point that the inclusion of an additional horizontal dimension poses no additional difficulties, while the procedure strongly relies on the linearity of the vertical structure. The first method that we discuss is a vertical average approach based on vertically constant levels, defining the interface interaction in a way that allows us to enforce some conservation-like properties. Secondly, we approximate using a $\mathbb{P}^1$ continuous finite element expansion, and in a familiar formulation we include layering via a pseudospectral collocation method. We finally remove the continuity assumption in the finite element setting and we perform layering with piecewise $\mathbb{P}^0$ and and $\mathbb{P}^1$ space-discontinuous Galerkin methods; note that in this latter case, the technique is applicable to nonlinear problems in a tractable way. For more details on the formulation and application of such techniques for the numerical approximation of time-dependent problems please refer to \cite{quarteroni} and \cite{donea} (and references therein).

\subsection{A ``conservative'' vertical averaging (CVA) approach }
We start by defining a set of uniformly distributed vertical levels $z_i,\, i=0\ldots N$, with $z_0=0$ and $z_N=H_z$ . For every interval $I_i=[z_{i-1},z_{i}]$ we define the vertically averaged quantity
\begin{equation}
q_i=\frac{1}{\dz} \int_{I_i}q(x,z,t)\, dz,
\end{equation}
where $\dz=H_z/N$.
We consider a local equation for every $I_i$ and we proceed by vertically integrating it, obtaining an exact formula for $q_i$:
\begin{equation}
\pt q_i + a\px q_i+\frac{1}{\dz}\,b\, \{\,q(z_i)-q(z_{i-1})\,\}+s\, q_i=0,
\end{equation}
We focus our attention on the approximation of the interface values $q(z_i)$ and $q(z_{i-1})$. We proceed by assuming a linear vertical variation between these two quantities within $I_i$ and imposing a consistency condition with the vertical average
\begin{equation}
\frac{1}{\dz}\int_{I_i}\frac{q(z_i)-q(z_{i-1})}{\dz}(z-z_i)+q(z_{i-1})\,dz=q_i,
\end{equation}
which after integration is reduced to
\begin{equation}
q(z_i)+q(z_{i-1})\approx 2q_i.
\end{equation}
This way to approximate interface values gives a coupling between the $N$ levels, and therefore the original 2D system is reduced to a one-dimensional system with a source term that mimics the vertical interaction; for instance, if we impose $q(z_0)=0$ in order to close the set of equations, this approach leads to
\begin{equation}
\pt q_k + a\px q_k+b\C q_k+s\, q_k=0,\qquad q_k=[\,q_0\ldots q_N]^t,
\end{equation}
with
\begin{equation}
\C_{i,j}=\frac{1}{\dz}
\left\{
\begin{array}{cc}
2 & \mbox{if } i=j,\\
(-1)^{i+j}\,4 & \mbox{if }  i>j,\\
0 & \mbox{if } i<j.
\end{array}
\right.
\end{equation}
Note that such formulation generates a source term $\C$ which, for a fixed mesh parameter $\dz$, has only positive eigenvalues $\lambda=\frac{2}{h}$.

\subsection{Piecewise linear finite element (FEM) expansion}
A different layered formulation can be obtained via a Galerkin ansatz in the vertical direction.
Proceeding in a similar manner as in the method of lines, we assume that $q(x,z,t)$ can be approximated by an expansion of basis functions in the vertical direction with space-time dependent coefficients,
\begin{equation}
q(x,z,t)\approx q_{\dz}= \sum_k q_k(x,t)\,\Phi_k(z).
\end{equation}
For a fixed pair $(x,t)$, we denote by $L^2(0;H_z)$ the space of square integrable functions in $z$ with its associated inner product
\begin{equation}\label{ip}
\lb u,v\rb_{L^2(0;H_z)}\equiv\int_0^{H_z} u(x,z,t)\,v(x,z,t)\,dz,
\end{equation}
and standard space $H^1(0;H_z)$ is defined upon the inner product
\begin{equation}
\lb u,v\rb_{H^1(0;H_z)}\equiv \lb u,v\rb_{L^2(0;H_z)}+ \lb \pz u, \pz v\rb_{L^2(0;H_z)}.
\end{equation}
We uniformly partition the interval $[0,H_z]$ into $N$ elements $I_i=[z_{i-1},z_i]\, i=1\ldots N$ and mesh parameter $\dz=H_z/N$. The approximation space $\mathcal{V}_{\dz}\subset H^1(0;H_z)$ is defined as
\begin{equation}
\mathcal{V}_{\dz}=\{v\in H^1(0;H_z):v|_{I_i}\in\mathbb{P}^1,\,i=1,\ldots,N\}.
\end{equation}
The set $\{\,\Phi_k(z)\,\}_{k=0}^N$ of hat functions constitutes a basis of the finite dimensional space $\mathcal{V}_{\dz}$. We write a weak formulation for (\ref{ley}),
\\
\textsl{Find $q_{\dz}\in\mathcal{V}_{\dz}$ such that}
\begin{equation}
\lb(\pt +a \px + s)\, q_{\dz},v\rb_{L^2(0;H_z)} + \lb b \pz q_{\dz},v\rb_{L^2(0;H_z)} =0\,,\qquad\forall v\in\mathcal{V}_{\dz}.
\end{equation}
In the same manner as in our first approach, this formulation is formally reduced to a 1D linear hyperbolic system for the expansion coefficients $\{q_k\}$:
\begin{equation}
\mathbf{M}\,(\pt q_k + a\px q_k +s q_k)+b\,\mathbf{K}\, q_k=0, \qquad q_k=[\,q_0\ldots q_ N\,]^t,
\end{equation}
where
\begin{equation}\label{matricesfem}
\mathbf{M}_{i,j}=\lb\Phi_i,\Phi_j\rb_{L^2(0;H_z)},\qquad \mathbf{K}_{i,j}=\lb\pz\Phi_i,\Phi_j\rb_{L^2(0;H_z)}.
\end{equation}
Since the inverse of $\mathbf{M}$ has to be computed only once, the simplified version reads
\begin{equation}
\pt q_k + a\px q_k + b \C q_k + \,s\, q_k=0,\qquad \C=\mathbf{M}^{-1}\mathbf{K}.
\end{equation}

Although the two layering approaches so far presented lead to one-dimensional hyperbolic systems of balance laws with a similar structure, the representation of an approximate two-dimensional solution is radically different. The first method will allow to recover piecewise constant average values for each vertical level, while the finite element method, with this particular choice of basis functions will be able to recover a piecewise linear and continuous approximate solution of the two-dimensional problem. In the same spirit, the forthcoming section develops an expansion in the $z$ direction such that a polynomial of arbitrary degree can be recovered.

\subsection{Pseudospectral collocation (PS) method}
The formulation of a pseudospectral method can be interpreted as a Galerkin approximation when a weight in the inner product is introduced and specific quadrature points are chosen. We begin by mapping the domain $\Omega$ into $\Omega_p=[0,L_x]\times[-1,1]$, where we consider again an expansion of the form $q_{p}=\sum_k\, q_k(x,t)\,\Phi_k(z)$, with a basis given by the set of Tschebychev polynomials $\{T_k(z)\}$ up to degree $p$; it is well-known that such polynomials are an orthogonal set in under the weighted inner product given by
\begin{equation}
\langle u,v\rangle =\int_{-1}^1\frac{1}{\sqrt{1-z^2}}\,u(x,z,t)\,v(x,z,t)\,dz.
\end{equation}
Given $p$, we define $p+1$ collocation points $z_i$ which in our case are optimally chosen as the zeroes of the Tschebychev polynomial $T_{p+1}(z)$, namely
\begin{equation}
z_i=-\cos\left[\frac{(2i+1)\pi}{2(p+1)}\right]\,,\qquad i=0,\ldots,p.
\end{equation}
We obtain a set of $p+1$ equations by imposing an \textsl{exact} solution at every collocation point, i.e.,
\begin{equation}
(\pt+a\px+b\pz+s)\sum_{k=0}^p\,q_k(x,t)\,\Phi(z)=0\,,\qquad \mbox{at every}\quad  z=z_i.
\end{equation}
The reduced system now is given by
\begin{equation}
\pt q_k+a\px q_k +b \C q_k+s\, q_k=0,\qquad q_k=[\,q_0\ldots q_ p\,]^t,
\end{equation}
with
\begin{equation}
\mathbf{M}_{i,j}=T_j(z_i), \qquad \mathbf{K}_{i,j}=\pz T_j\,(z_i), \quad\mbox{and}\quad \C=\mathbf{M}^{-1}\mathbf{K}.
\end{equation}

\subsection{Discontinuous Galerkin (DGP0 - DGP1) approximation}
Proceeding in a similar way as in the piecewise continuous and linear finite element approach, after partitioning the vertical domain we consider a local problem in the interval $I_i=[z_{i-1},z_{i}]$, which is mapped into the reference element $I_{ref}=[-1,1]$ via
\begin{equation}
z=\frac{z_{i-1}+z_{i}}{2}+\xi\,\frac{\dz}{2}.
\end{equation}
In such domain we consider a local expansion
\begin{equation}
q_i(x,z,t)\approx \sum_j c_{j,i}(x,t)\,\Phi_j(z),
\end{equation}
where the basis $\{\,\Phi_j(z)\,\}_{j=0}^p$ is chosen as the set of Legendre polynomials up to degree $p$ which is orthogonal over $[-1,1]$. If $p=1$, then $q_i(x,z,t)=c_{0,i}+\xi\, c_{1,i}$ and after multiplication by test functions in the space of polynomials of degree 1 over $[-1,1]$, $\mathbb{P}^1_{I_{ref}}$, and integration by parts, we write the following weak formulation (with the same $L^2$ product as in (\ref{ip}) but in the interval $[-1;1]\,$) for the local problem :
\textsl{Find $q_i\in\mathbb{P}^1_{I_{ref}}$ such that}
\begin{equation}
\lb (\pt+a \px+s)\, q_i,v\rb_{L^2(-1;1)} +v\, b \,q\, |_{-1}^1-\lb b\, q_i,\pz v\rb_{L^2(-1;1)} =0\qquad\forall v\in\mathbb{P}^1_{I_{ref}}.
\end{equation}
Note that the basis in the expansion are no longer of compact support as in the piecewise linear continuous finite element approach, and therefore there is no continuity enforcement; every local problem is treated separately. However, we do enforce coupling between the different layers by  a suitable approximation of the interface values arising in the term $v b q_i |_{-1}^1$; in this article we approximate this interface flux via upwind considerations. Before focusing our attention on this issue, we recall that testing is done with $v=1$ and $v=\xi$, leading to the following set of equations for the coefficients $c_{0,i}$ and $c_{1,i}$ at every element:
\begin{eqnarray}
\pt c_{0,i}+a\px c_{0,i}+s\, c_{0,i}+\frac{1}{\dz}\,\{\,b\,q(z=z_i)-b\,q(z=z_{i-1})\,\}&=&0,\\
\pt c_{1,i}+a\px c_{1,i}+s\, c_{1,i}+\frac{3}{\dz}\,\{\,b\,q(z=z_i)+b\,q(z=z_{i-1})-2b\,c_{0,i}\,\}&=&0.
\end{eqnarray}
This system of equations can be closed by suitable boundary conditions on vertical extremes of the domain; the remaining issue is henceforth the approximation of the interfaces values. We use an upwind flux, that for $b>0$ leads to:
\begin{eqnarray}
b\,q(z=z_i)&=&b\,q_{i}(\xi=1)=b\,( c_{0,i}+c_{1,i} )\\
b\,q(z=z_{i-1})&=& b\,q_{i-1}(\xi=1)=b\,( c_{0,i-1}+c_{1,i-1}).
\end{eqnarray}
As in all our previous approaches, the dimensionally reduced set of equations is a one-dimensional system of hyperbolic balance laws, consisting in this case in two equations determining the local coefficients $c_{0,i}$ and $c_{1,i}$ for every element; there is an increasing computational cost depending on the local approximation degree compared with continuous formulations.

The next section takes as a starting point the resolution of a hyperbolic system of balance laws in the form of those obtained by the layering procedures so far introduced; they all share a single horizontal advection speed while a stable reactive source establishes coupling and transmission between the different vertical levels. The same approach can be extended to systems where different horizontal advection speeds are present, as we will show in the examples, but since it is possible to work in the characteristic variable field they key points remain in the resolution of this simpler problem.

\section{The ADER approach for high horizontal resolution}
This section consists of a concise version of the original scheme formulation; it makes extensive use of the linearity of the system while the general formulation has been developed by conservation laws in a wide framework. For more details about the general formulation please refer to \cite{ader-1,ader-2,toro} and references therein.
We consider a system of balance laws in one spatial dimension given by
\begin{equation}\label{sbl}
\pt q +\A \px q+\S q=0,
\end{equation}
along with initial and boundary conditions. We will assume $\A$ diagonal, and $\S$ a definite semipositive matrix. We partition the horizontal domain into intervals $E_j=[x_{j-1/2},x_{j+1/2}]$, and  we consider  a control volume in the $x-t$ space $E_j\times [t^{n},t^{n+1}]$  with dimensions $\dx=x_{j+1/2}-x_{j-1/2}$, $\dt=t^{n+1}-t^n$. Integrating(\ref{sbl}) over the control volume leads to
\begin{equation}\label{abl}
q_j^{n+1}=q_j^n-\frac{\dt}{\dx}\,\left(f_{j+1/2}-f_{j-1/2}\right)\,-\dt \,s_j,
\end{equation}
where
\begin{equation}\label{fluxsource}
f_{+1/2}=\frac{1}{\dt}\int_{0}^{\dt}\!\!\!\!\A\, q(x_{j+1/2},\tau)\,d\tau,\quad s_j=\,\frac{1}{\dt}\frac{1}{\dx}\int_{0}^{\dt}\!\!\!\!\int_{x_{j-1/2}}^{x_{j+1/2}}\!\!\!\! \S q(x,\tau)\,dx\,d\tau ,
\end{equation}
which provides a marching exact formula for the spatially averaged quantity
\begin{equation}
q_j^n=\frac{1}{\dx}\int_{x_{j-1/2}}^{x_{j+1/2}}\!\!\!\!q(x,t^n)\,dx.
\end{equation}
The ADER scheme is a numerical method that generates a high-order approximation of the two integrals appearing in (\ref{fluxsource}) and consequently, a high-order approximation of the averaged state value. The key ingredients of this scheme applied to this system are: 1) a reconstruction procedure that allows us to recover high-order degree polynomials at every element given the set of averaged values, for 2) the formulation of a generalized Riemann problem for, 3) the generation of Taylor series expansions for the approximation of the required states. After integration is performed, the resulting scheme is of order $m$, both in space and time, where $m-1$ corresponds to the polynomial degree used in the reconstruction step, provided that the integration procedures are of at least the same order of accuracy; in our case, time integrals will be treated exactly, while the space integral remaining in the source term will be approximated by a five point Gaussian quadrature (in this article we study schemes up to order 4).
For the reconstruction step, we the use WENO procedure proposed in \cite{dumbweno}, which slightly differs from the original version \cite{osher,shu1}, as it recovers a single polynomial for every cell, which can be used either for extrapolation at the boundary or evaluation inside the cell. As we will be dealing with systems of balance laws, this latter characteristic is useful for the approximation of the source term. Henceforth we denote by ADERm to the scheme that uses a $m-1$ polynomial degree reconstruction.
We now focus on the task of approximation of the fluxes in (\ref{abl}), which requires the approximation of the state $q$ at the interface points in every horizontal element, i.e, $q(x_{j+1/2},\tau)$ and $q(x_{j-1/2},\tau)$. After every time step inicialization at $t^n$, being performed the reconstruction step, for every element we dispose of a reconstruction polynomial $p_j(x)$ of degree $m-1$. At every cell interface we pose the following generalized Riemann problem:
\begin{eqnarray}\label{grp}
PDE:\;&  \pt q +\A\px q+\S q=0,\\
IC:\; & q(x,0)=\left\lbrace\begin{array}{cc}q_L(x)=p_j(x), & x<x_{j+1/2}\\q_R(x)=p_{j+1}(x), & x>x_{j+1/2}\end{array}\right.
\end{eqnarray}
The ADER approach links the resolution of this generalized Riemann problem with the approximation of the interface state $q(x_{j+1/2},\tau)$. Assuming $\tau$ sufficiently small, we write a Taylor expansion in time for the interface state
\begin{equation}
q(x_{j+1/2},\tau)\approx q(x_{j+1/2},0^+)+\sum_{k=1}^{m-1}[\pt^{(k)}q(x_{j+1/2},0^+)]\frac{\tau^k}{k!},
\end{equation}
where the time derivatives are replaced with space derivatives via the Cauchy-Kowalevski procedure ( also know as Lax-Wendroff procedure \cite{lax}) procedure. In the simplest case, if $\A$ is a diagonal matrix with a single element, the $k-th$ time derivative holds the binomial formula
\begin{align}\label{bin}
\pt^{(k)}q(x,t)&=(-1)^k(\A\px+\S)^kq(x,t)\,,\\
&=(-1)^k\sum_{i=0}^k\binom{k}{i}\A^i\S^{k-i}\px^{(k)}q\,,\qquad k=0,\ldots ,m-1.
\end{align}
The interface state and its spatial derivatives are approximated by the solution at $(x-x_{j+1/2})/\tau$=0 of the following set of classical Riemann problems derived from (\ref{grp}):
\begin{align}\label{crp}
PDE:\;&  \pt \px^{(k)}q +\A \px^{(k)}q=0,\qquad k=0,\ldots,m-1\\
IC:\; & \px^{(k)}q(x,0)=\left\lbrace\begin{array}{cc}\px^{(k)}q_L(x_{j+1/2})=\px^{(k)} p_j(x_{j+1/2}), & x<x_{j+1/2}\\\px^{(k)}q_R(x_{j+1/2})=\px^{(k)} p_{j+1}(x_{j+1/2}), & x>x_{j+1/2}\end{array}\right..
\end{align}
There are two things worth to point out at this stage: First, that given $\A$ a diagonal matrix with a single entry, the formula (\ref{bin}) uses the fact that $\A$ and $\S$ commute; the case with multiple entries in the diagonal is still tractable even though it is not possible to make use of the binomial formula; this same diagonal structure allows an exact solution of (\ref{crp}) component-wise where it is reduced to a linear advection problem, while the more general formulation of the ADER scheme include the use of approximate Riemann solvers at this point. Secondly, we note that although the set of classical Riemann problems used to approximate $q(x_{j+1/2},0^+)$ and its space derivatives neglect the presence of the source term, it is still captured in the final Taylor expansion

\begin{equation}\label{tf}
q(x_{j+1/2},\tau)\approx q(x_{j+1/2},0^+)+\sum_{k=1}^{m-1}\left\{\sum_{i=0}^k\binom{k}{i}\A^i\S^{k-i}\px^{(k)}q(x_{j+1/2},0^+)\right\}\frac{(-\tau)^k}{k!}.
\end{equation}

We insert this formula into the flux (\ref{fluxsource}) and integrate, yielding
\begin{equation}
f_{j+1/2}\approx \A\left( q(x_{j+1/2},0^+)+\sum_{k=1}^{m-1}\left\lbrace\sum_{i=0}^k\binom{k}{i}\A^i\S^{k-i}\px^{(k)}q(x_{j+1/2},0^+)\right\rbrace\frac{(-\dt)^k}{k+1!}\right).
\end{equation}

The remaining issue on the scheme is the approximation of the source term integral $s_j$; as we have previously seen, since we write Taylor expansions in time and the source is linear, we don't have to worry about time integration as it can be done exactly; the space integral is approximated by a Gaussian quadrature rule, therefore
\begin{equation}
s_j=\,\frac{1}{\dt}\frac{1}{\dx}\int_{0}^{\dt}\!\!\!\!\int_{x_{j-1/2}}^{x_{j+1/2}}\!\!\!\! \S q(x,\tau)\,dx\,d\tau\approx
\,\frac12\frac{1}{\dt}\int_{0}^{\dt}\sum_{\alpha=1}^{N_{g}} w_\alpha\,\S q(x_{\alpha},\tau)\,d\tau
\end{equation}
where $N_{g}$ stands for the number of Gauss points $\{x_{\alpha}\}_{\alpha=1}^{N_{g}}$ used in the integration rule; its associated weights are denoted by $w_{\alpha}$. Within the $j-th$ element, we write a Taylor series in time and replace the time derivatives by space derivatives in a similar way as it has been done for the calculation of the fluxes; therefore, inside the element $E_j$, in a given Gauss point $x_{\alpha}$ the Taylor formula (\ref{tf}) holds replacing $x_{j+1/2}$ by $x_{\alpha}$. What is different however, is that the state values $\px^{(k)}q(x_{\alpha},0^+)$ are not obtained from the solution of a Riemann problem, but from the solution of the Cauchy problem at the $j$-th cell
\begin{align}
PDE:\;&  \pt \px^{(k)}q +\A \px^{(k)}q=0,\qquad k=0,\ldots,m-1\\
IC:\; & \px^{(k)}q(x,0)=\px^{(k)}p_j(x).
\end{align}

Performing the whole procedure (reconstruction, flux and source term calculation, and update) leads to a scheme of $m-th$ order of accuracy both in space and time, with an optimal CFL condition of $CFL\leq 1$ without source term; in certain cases, as we will see, the source term does not affect the stability condition and can consequently be omitted when determining the time step, while in other situations its incorporation restricts the time stepping in a way that other solutions, as splitting schemes are required; in any case, the ADER method provides an accurate and efficient way to treat the advective character of the studied equations.

We now perform numerical tests for the global scheme, i.e., the numerical procedure that includes a low-order layering in the vertical direction generating a system of balance laws that is approximated by an ADER$m$ scheme. The purpose of the tests are twofold: one on hand they validate the accuracy of the proposed ADER$m$ scheme in the presence of different source terms, while they also study the performance of the global scheme and its properties.

\section{Numerical experiments}
In this section we consider two models in order to test the scheme presented in this article. For both cases we consider  a two-dimensional $x-z$ domain $\Omega=\{(x,z)\in[0,L]\times[0,H]\}$. In every case numerical convergence rates are shown in order to assess the accuracy of the scheme, while other relevant tests dealing with stability, conservation and physical validation are performed. Most of the code concerning the reconstruction procedure, as well as flux and source term calculations were coded and vectorized in \verb"FORTRAN95"; simulations were performed using an Intel Core 2 Duo 2.50 Ghz processor in a computer with 4 Gb RAM.

\subsection {2D linear advection}
The simplest case in which our approach is meaningful consists on the advection of a quantity $q=q(x,z,t)$ through $\Omega$. This phenomenon is represented by the linear, two dimensional advection equation
\begin{equation}
\pt q+U\px q+ W\pz q=0,
\end{equation}
where the pair $(U,W)$ is a constant two-dimensional velocity field. The initial condition $q(x,z,0)$ is denoted by $q_0(x,z)$. Well-posedness of this problem requires that for $U$ and $W$ positive, q is specified at $x=0$ and $z=0$; we set $q(0,z,t)=q(L,z,t)$ and $q(x,0,t)=q_0(x,0)$.
Given such model, we first perform vertical layering via the  methods stated in section 2; as it was previously shown, every case leads to a system
\begin{equation}
\pt q_k+U\px q_k +\C q_k=0,
\end{equation}
where the resulting source term is obtained by replacing $b =W$. This hyperbolic system is solved with the ADER approach described in section 3.

\subsubsection{Convergence studies}
We first perfom numerical tests in order to validate the convergence rates of the ADER method for a system in the presence of a source term. Every source term in this case is stable in the sense that its eigenvalues are contained in the left half-plane. We set $L_x=H_z=1$, $U=W=0.01$ and an initial condition \begin{equation}q_0(x,z)=sin(10\pi x)\,sin(2\pi z).\end{equation} Time stepping  is established as $\dt=CFL\,(\dx/U)$, with a CFL number chosen to be 0.8, being the simulations stable with CFL up to 1; note that the time stepping is computed  as for an homogeneous linear advection problem in one dimension. Final time for simulation is $t=50[s]$. For 10 vertical levels (or collocation points in the case of the PS source term), and with ADER order varying from 2 to 4 we obtain the results shown in Tables \ref{t1}-\ref{t2}, where relative errors are computed in space at final time of simulation with respect to a reference solution obtained with 1024 horizontal elements. Predicted convergence rates for smooth solutions, in the presence of every source term proposed in section 2, are achieved in both $L_1$ and $L_{\infty}$ norms.

\begin{table}
\centering
\caption{Linear advection: convergence rates and relative error ($\%$) for the ADER$m-th$ horizontal order scheme with a conservative vertically averaged (ADER$m$-CVA) and piecewise constant discontinuous Garlerkin source terms (ADER$m$-DGP0).}\label{t1}\vskip 8mm
\begin{tabular}{l c c c c c}
\hline\\
Method & N & $L_{\infty}$ error &  $L_{\infty}$ order & $L_1$ error &  $L_1$ order \\ [0.5ex]
\hline\\ \\
ADER2-CVA& 32 &  3.3429e-002 &  &7.5757e-003& \\
& 64 &   1.2924e-002 & 1.3711 & 2.9740e-003&1.9893 \\
& 128 &  5.0697e-003  & 1.3500 &7.5953e-004& 1.9692\\
& 256 &   2.0101e-003 &1.3347  &1.9504e-004& 1.9614  \\
& 512 & 8.0148e-004 & 1.3265 &5.0518e-005 & 1.9489 \\ \\

ADER3-CVA & 32 & 8.2992e-004 &  &   1.1808e-002& \\
& 64 &   9.7975e-005& 3.0825 & 6.2395e-005 &3.0906  \\
& 128 &  1.1854e-005 & 3.0471 &  7.5473e-006&  3.0474\\
& 256 &  1.4568e-006   &3.0244  & 9.2749e-007 & 3.0246 \\
& 512 &  1.8056e-007 & 3.0123 &1.1495e-007& 3.0123 \\ \\

ADER4-CVA & 32 &  4.6056e-005 &  &3.0473e-005& \\
& 64 &   2.9604e-006 & 3.9595 & 1.8899e-006&4.0112 \\
& 128 &  1.8480e-007& 4.0018 & 1.1738e-007& 4.0091\\
& 256 &  1.1457e-008 &4.0117  & 7.2826e-009& 4.0106\\
& 512 &  6.7165e-010 & 4.0923&4.2755e-010& 4.0903\\\\ \\

ADER2-DGP0 & 32 & 2.3473e-002 &  &5.9531e-003& \\
& 64 &   9.7303e-003&  1.9893  &1.5895e-003&1.9051 \\
& 128 &  3.9960e-003  &1.9692  &4.2553e-004&   1.9012\\
& 256 &   1.6320e-003 &1.9614    &  1.1495e-004& 1.8883   \\
& 512 &  6.6328e-004 & 1.9489 &3.0676e-005& 1.9058\\ \\

ADER3-DGP0 & 32 & 4.2071e-004 &  &  2.6867e-004& \\
& 64 &   5.2219e-005& 3.0906 &3.3291e-005&3.0126  \\
& 128 &  6.4997e-006 & 3.0474 &  4.1373e-006& 3.0084\\
& 256 &  8.1006e-007  &3.0246  & 5.1572e-007 & 3.0040 \\
& 512 &  1.0111e-007 & 3.0123 & 6.4369e-008& 3.0022 \\ \\

ADER4-DGP0 & 32 &  8.2443e-006 &  &4.9867e-006& \\
& 64 &5.7959e-007& 4.0112 & 3.6692e-007&3.7646 \\
& 128 &  3.8606e-008& 4.0091 & 2.4559e-008& 3.9011\\
& 256 & 2.4817e-009 &4.0106  & 1.5797e-009& 3.9586\\
& 512 &  1.4823e-010 & 4.0903&9.4370e-011& 4.0652\\

\end{tabular}
\end{table}

\begin{table}
\centering
\caption{Linear advection: convergence rates and relative error ($\%$) for the ADER$m-th$ horizontal order scheme with a piecewise linear finite element (ADER$m$-FEM) and pseudospectral source terms (ADER$m$-PS).}\label{t2}\vskip 8mm
\begin{tabular}{l c c c c c}
\hline\\
Method & N & $L_{\infty}$ error &  $L_{\infty}$ order & $L_1$ error &  $L_1$ order \\ [0.5ex]
\hline\\ \\

ADER2-FEM & 32 &  2.6751e-002 &  &7.5757e-003& \\
& 64 &   1.1372e-002& 1.2341 & 1.9927e-003&1.9267 \\
& 128 &  4.7237e-003  & 1.2675 & 5.3943e-004& 1.8852\\
& 256 &   1.9398e-003 &1.2840   &  1.4590e-004& 1.8865  \\
& 512 &  7.9072e-004 & 1.2947 &3.8829e-005& 1.9097\\ \\

ADER3-FEM & 32 & 7.2227e-004 &  &  4.6239e-004& \\
& 64 &   9.0743e-005& 2.9927 & 5.7778e-005&3.0005  \\
& 128 &  1.1330e-005 & 3.0017 &  7.2156e-006& 3.0013\\
& 256 &  1.4153e-006  &3.0009  & 9.0102e-007 & 3.0015 \\
& 512 &  1.7682e-007 & 3.0008 &1.1257e-007& 3.0007\\ \\

ADER4-FEM & 32 &  2.2579e-005 &  &1.4252e-005& \\
& 64 &   1.4553e-006 & 3.9556 & 9.2260e-007&3.9493 \\
& 128 &  9.2783e-008 & 3.9713 & 5.9028e-008& 3.9662\\
& 256 &  5.8391e-009 &3.9900  & 3.7170e-009& 3.9892\\
& 512 &  3.4532e-010 & 4.0798&2.1984e-010& 4.0796\\\\ \\

ADER2-PS & 32 &  9.3519e-003 &  &4.6919e-003& \\
& 64 &   2.6984e-003& 1.7931 & 1.1743e-003&1.9983 \\
& 128 &  9.6431e-004  & 1.4846  & 2.9331e-004 & 2.0014\\
& 256 &  3.6181e-004 &1.4143   &  7.3797e-005& 1.9908  \\
& 512 &  1.3932e-004 & 1.3768 &1.8711e-005& 1.9796\\ \\

ADER3-PS & 32 & 1.9786e-004 &  &  1.2654e-004& \\
& 64 &   2.3921e-005 & 3.0481 & 1.5244e-005 &3.0532  \\
& 128 &  2.9354e-006 & 3.0267 &  1.8692e-006& 3.0278\\
& 256 &  3.6341e-007  &3.0139  & 2.3137e-007 & 3.0141 \\
& 512 &  4.5205e-008 &  3.0071  &2.8779e-008& 3.0071\\ \\

ADER4-PS & 32 & 1.6061e-005 &  &1.0191e-005& \\
& 64 &  9.3375e-007 & 4.1044 & 5.9498e-007&4.0983 \\
& 128 & 5.6164e-008 &4.0553 & 3.5750e-008 & 4.0568\\
& 256 & 3.4266e-009 &4.0348  &2.1815e-009& 4.0346\\
& 512 & 1.9952e-010 & 4.1022&1.2702e-010& 4.1022\\

\end{tabular}
\end{table}

\subsubsection{Advection of a gaussian initial condition and averaging}

\begin{figure}
\caption{Linear advection of Gaussian profile: temporal evolution of the maximum value of the approximated solution for different source terms and horizontal orders. Maximum value is expected to exit the domain at $t=50[s]$.}
\begin{minipage}[b]{0.5\linewidth}
\centering
\includegraphics[height=0.25\textheight,width=\textwidth]{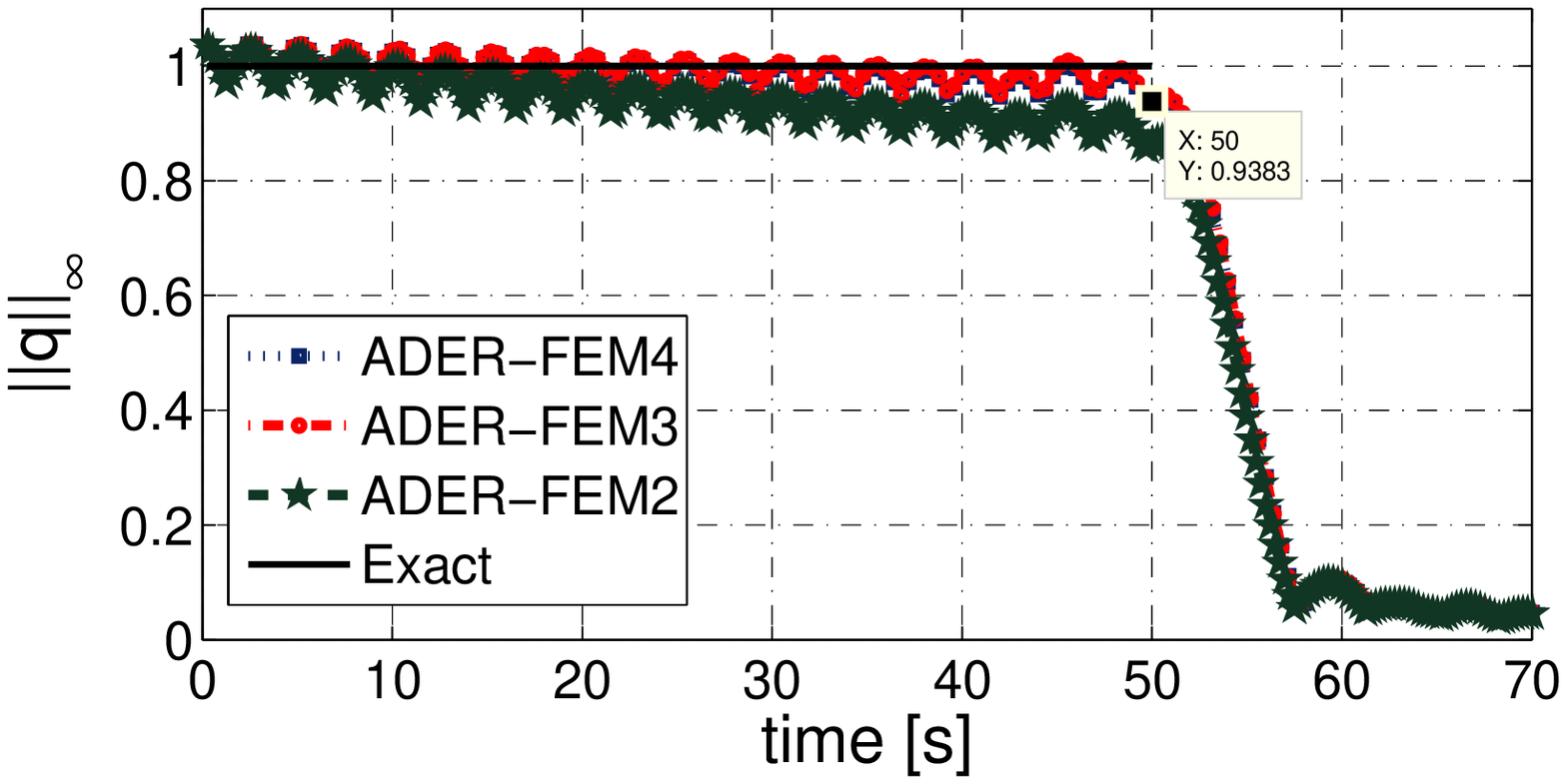}
\includegraphics[height=0.25\textheight,width=\textwidth]{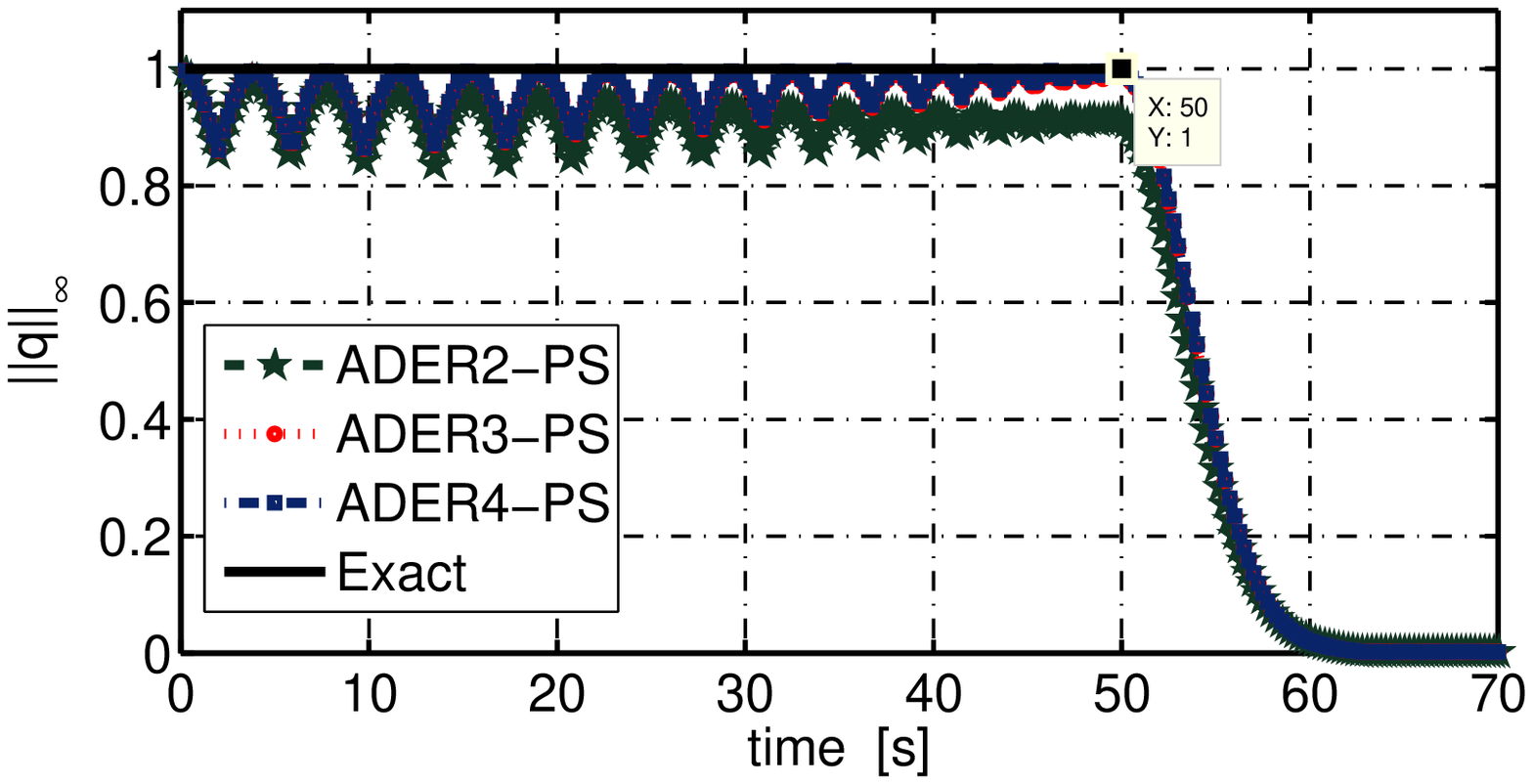}
\end{minipage}
\hspace{0.0cm}
\begin{minipage}[b]{0.5\linewidth}
\centering
\includegraphics[height=0.25\textheight,width=\textwidth]{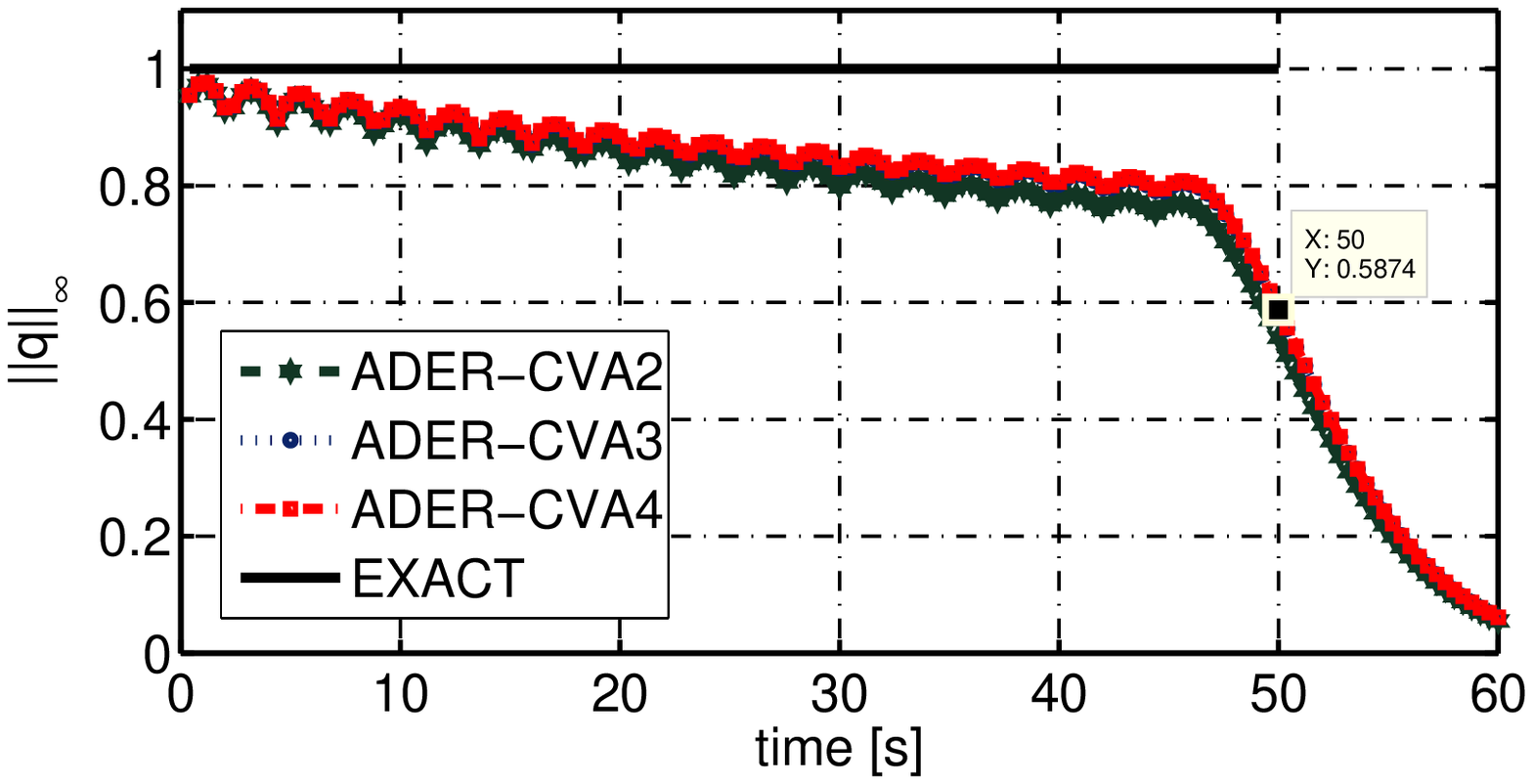}
\includegraphics[height=0.25\textheight,width=\textwidth]{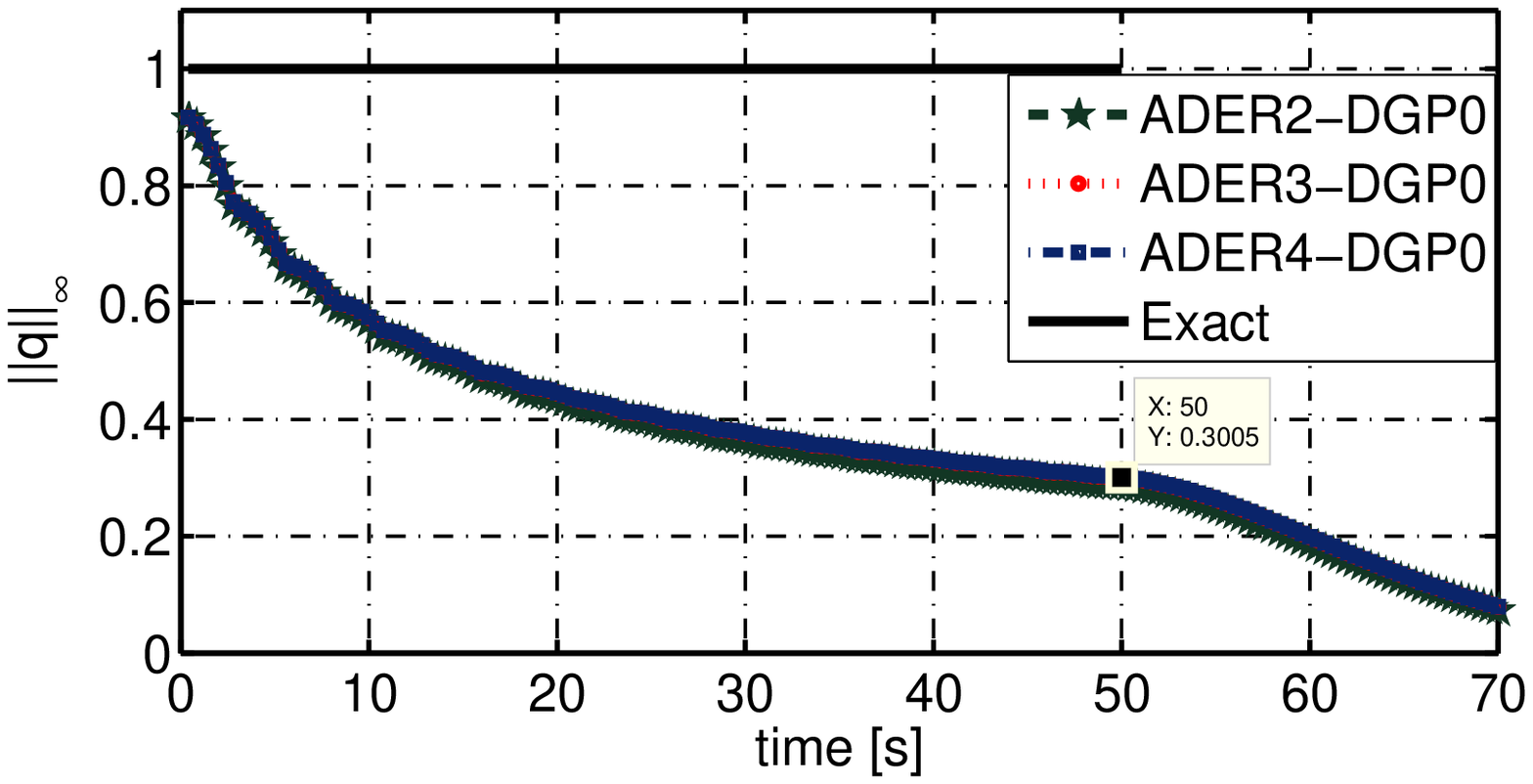}
\end{minipage}
\label{fig:linf1}
\end{figure}

In a second numerical experiment we consider the same spatial domain and velocity field, but an initial condition
\begin{equation} q_0(x,z)=\exp\left(-\frac{(x-0.5)^2+(z-0.5)^2}{0.05^2}\right),
\end{equation}
and a simulation time of $70[s]$ (with $100 \times 100$ elements). We measure the $L_{\infty}$ norm of the 2D solution in order to discuss the diffusive properties of the global scheme. We observe that the application of the ADER method in the horizontal is not causing major difussion when orders 3 or 4 are used, while the choice of the source term is relevant. Figure \ref{fig:linf1} shows that both FEM and PS source terms lead to schemes that are not diffusive at all while exhibiting an oscillation due to the transmission of the maximum value from one layer to another; the CVA approach shows considerable diffusion while the DGP0 scheme with upwind flux is the most diffusive one, wich was expected since it is the only first order method. Another aspect to be observed from these graphs is the ability to capture the correct advection speed; it can be noted that the CVA approach overestimates the physical speed, while the remaining methods clearly show a decay at $t=50 [s]$ that corresponds to the expected time for the top of the advected initial condition to leave the physical domain. The oscillations observed afterwards in the FEM source term case are due to a spurious reflection in the upper boundary.

\begin{figure}
\caption{Linear advection and averaging: 1.Temporal evolution of the vertically averaged exact solution $q_s$ (left) 2. Temporal evolution of the full subdomain average (right).}
\begin{minipage}[b]{0.5\linewidth}
\centering
\includegraphics[height=0.25\textheight,width=\textwidth]{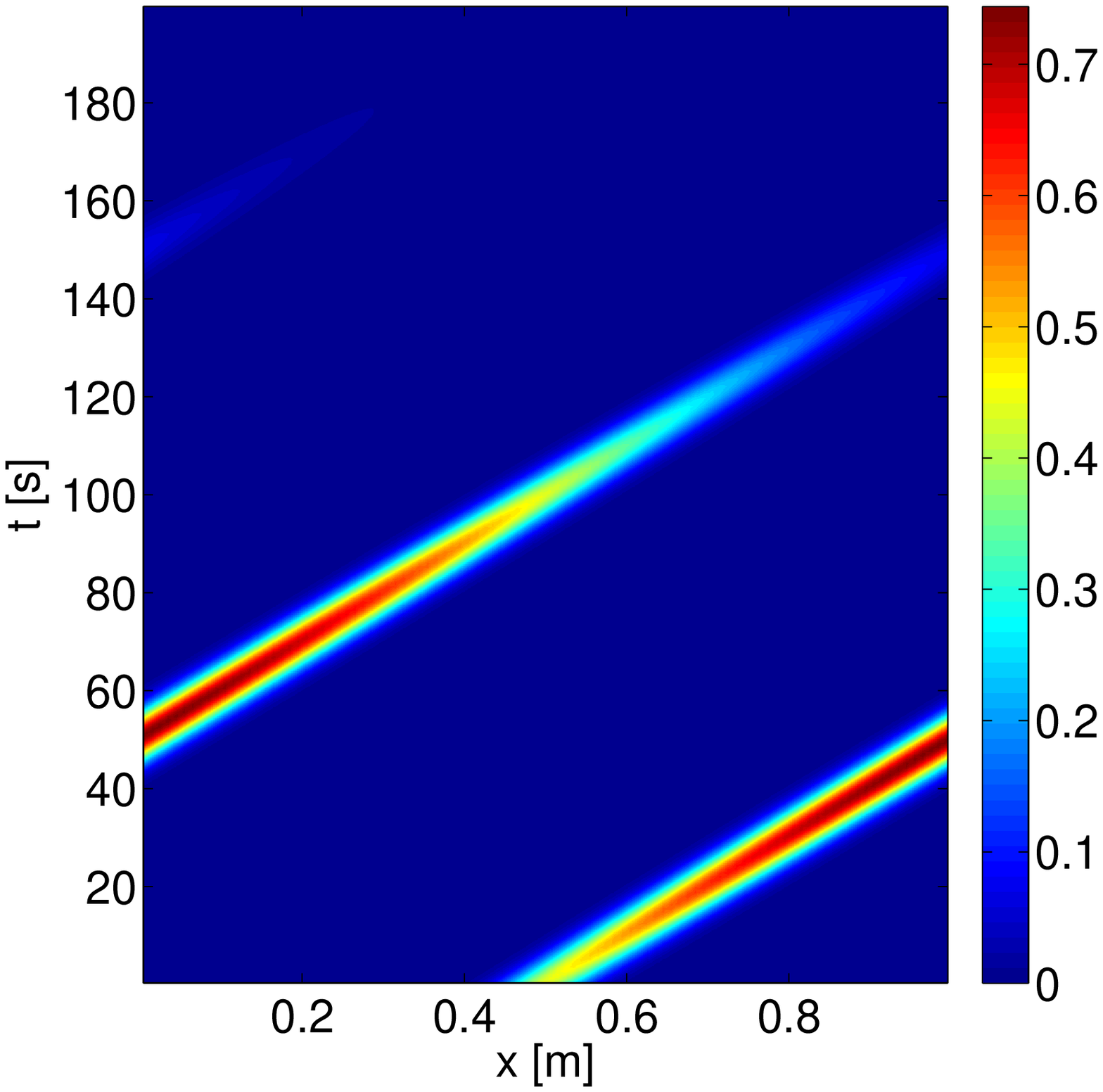}
\end{minipage}
\hspace{0.0cm}
\begin{minipage}[b]{0.5\linewidth}
\centering
\includegraphics[height=0.25\textheight,width=\textwidth]{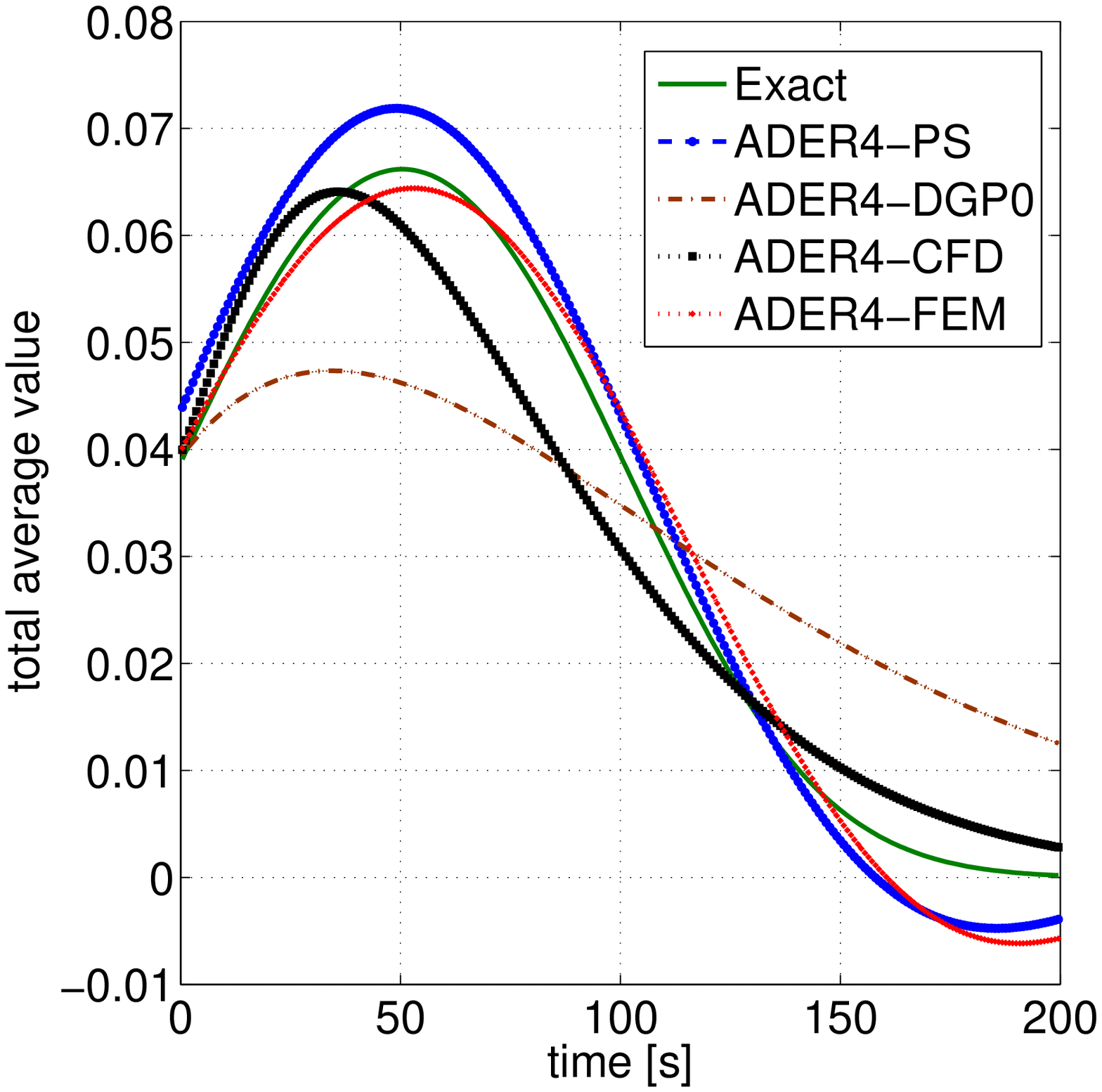}
\end{minipage}
\label{av1}
\end{figure}

\begin{figure}
\caption{Linear advection and averaging: vertically averaged solutions (right) and its residuals $|\,q_s^{exact}-q_s^{approx}\,|$ (left). $100\times 20$ elements. Top: ADER4-FEM. Bottom: ADER4-PS}\label{av2}
\centering
\vskip 5mm
\includegraphics[height=.3\textheight,width=\textwidth]{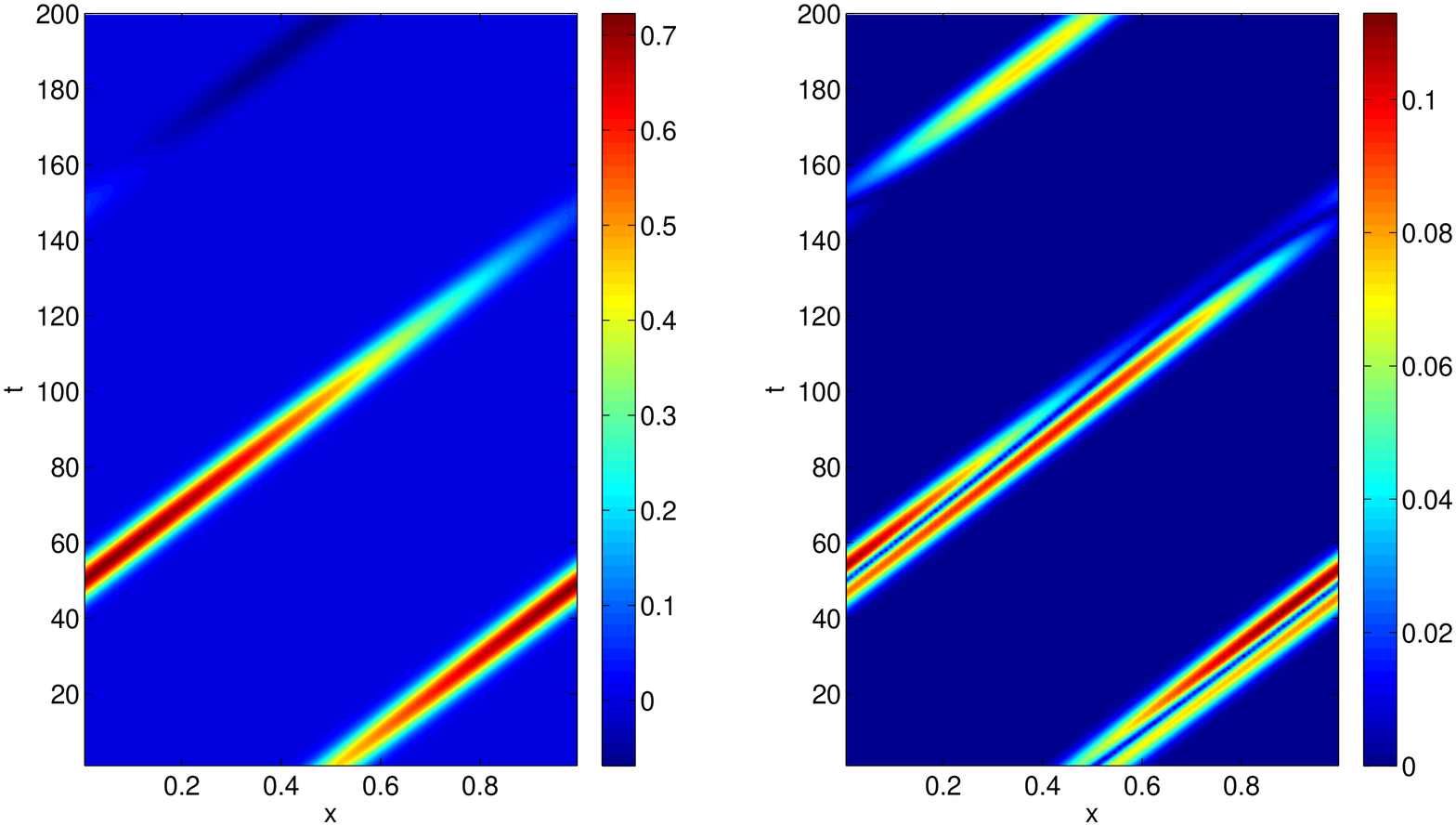}
\vskip 5mm
\includegraphics[height=.3\textheight,width=\textwidth]{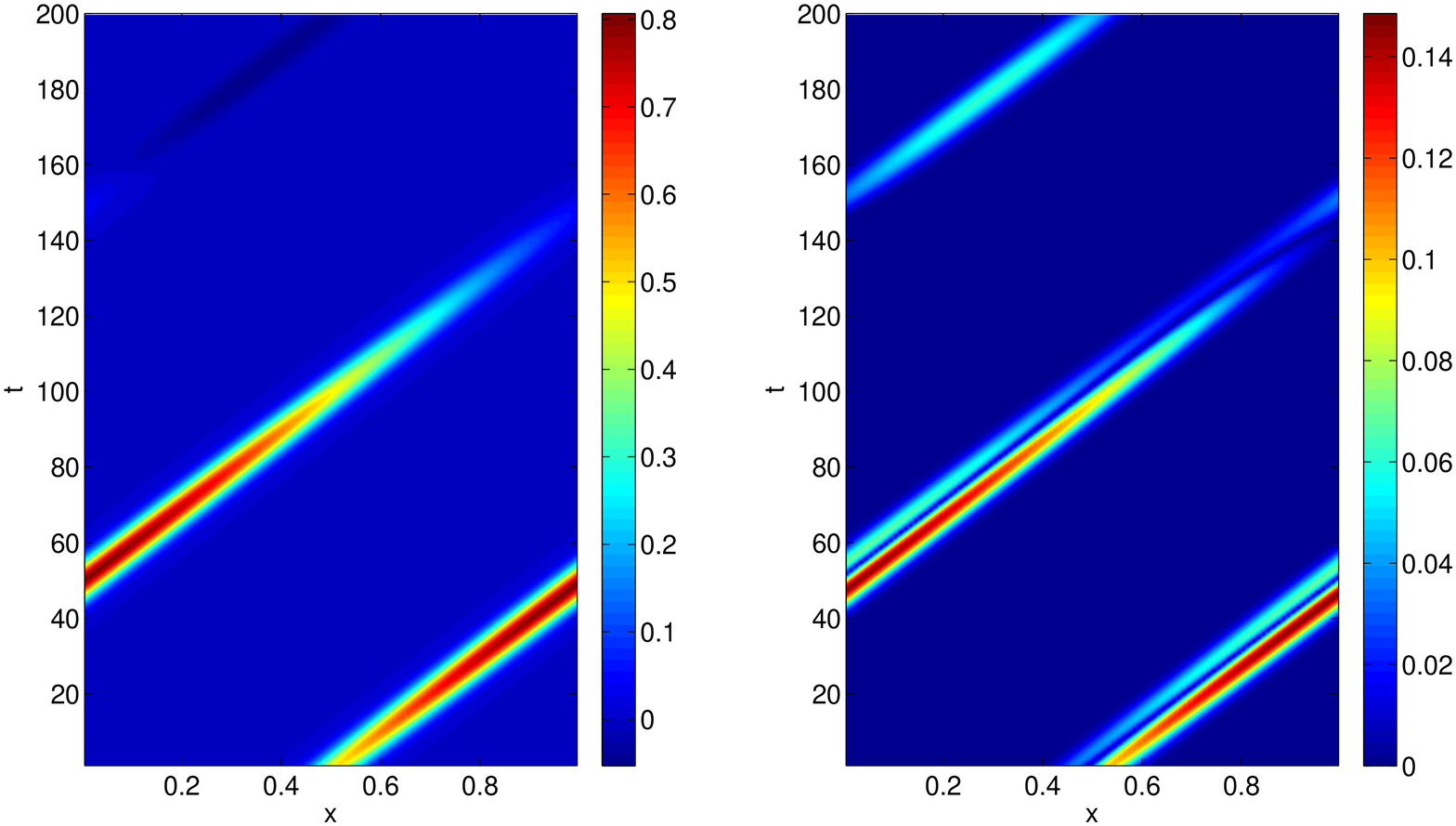}

\end{figure}

\begin{figure}
\caption{Linear advection and averaging: vertically averaged solutions (right) and its residuals $|\,q_s^{exact}-q_s^{approx}\,|$ (left). $100\times 20$ elements. Top: ADER4-CVA. Bottom: ADER4-DGP0}\label{av3}
\centering
\vskip 5mm
\includegraphics[height=.3\textheight,width=\textwidth]{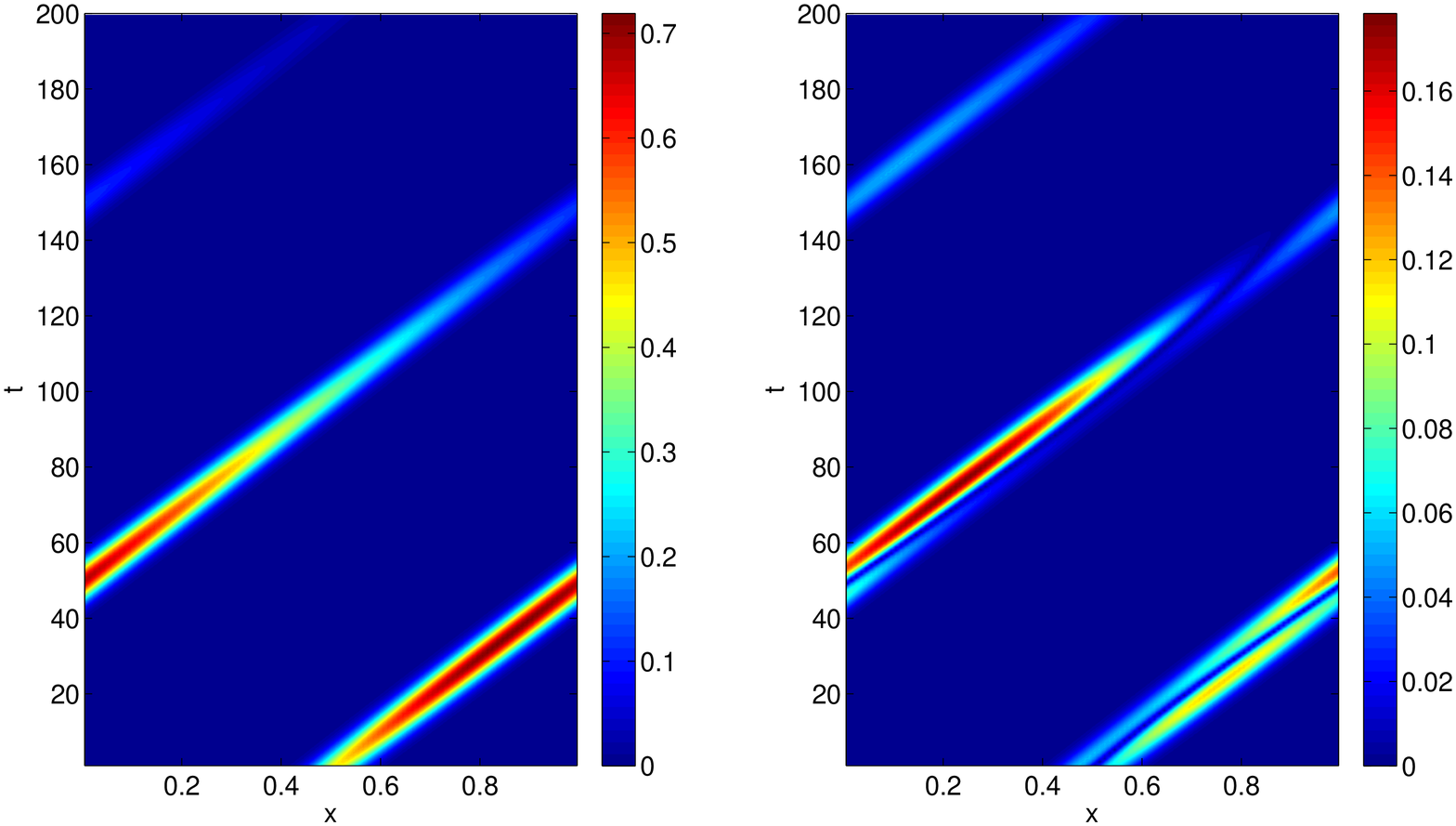}
\vskip 5mm
\includegraphics[height=.3\textheight,width=\textwidth]{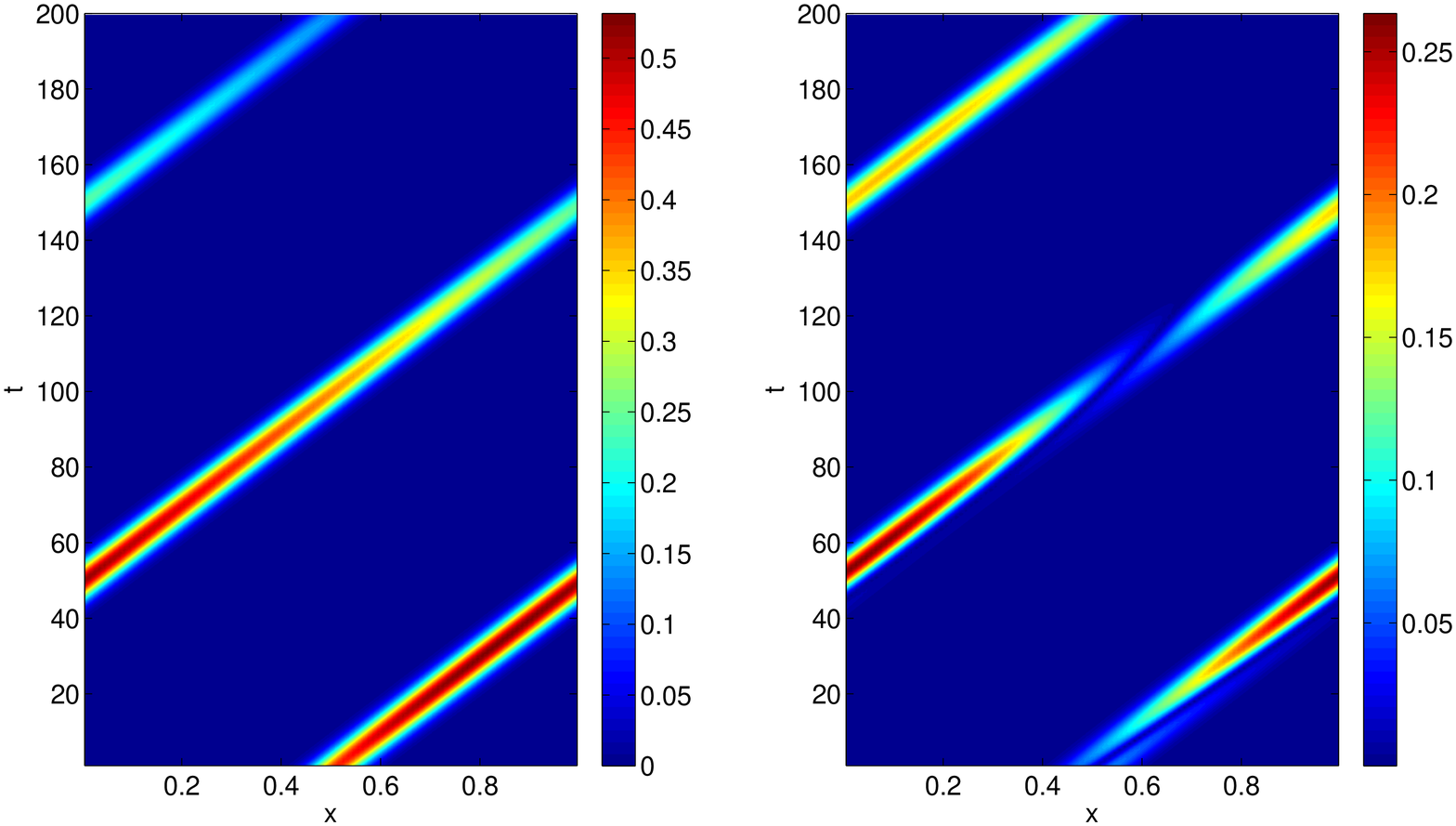}

\end{figure}

We conclude this example by analyzing a setting that is relevant in a data assimilation context. We consider the same Gaussian initial condition with $U=0.01$, $W=0.001$, but we also consider a subdomain $\Omega_s=[0,1]\times[0.5,0.6]$ where we vertically average the solution. This leads to a quantity
\begin{equation}
q_s=q_s (x,t)=\frac{1}{0.1}\int_{0.5}^{0.6}q(x,z,t)\,dz,
\end{equation}
which is computed upon the approximated solution. Simulation time is now $200[s]$ with 100 horizontal elements and 20 vertical elements. We show the evolution of $q_s$, which is advected at the same speed, and we show absolute residuals of the averaged values obtained from approximated solutions and the exact average value. Figure \ref{av1} shows the space-time evolution of the exact quantity $q_s$ and the temporal evolution of $q$ totally averaged over $\Omega_s$ with different layering approaches; it can be seen that again the CVA approach seems to overestimate the advection speed, the same situation happens with the DGP0 scheme, being the FEM source term slightly better than the PS on generating a scheme that tries to recover $q_s$; such performance is confirmed by Figures \ref{av2}-\ref{av3}, where the ADER4-FEM scheme is the one exhibiting the lower residuals.

\subsection{2D atmospheric model}

We also study the performance of our approach in a linear, non-hydrostatic atmospheric model. If we start from conservation principles, the model is usually expressed in terms of the velocity field in the $x-z$ plane $(u,w)$, the density $\rho$, the pressure $p$ and the temperature $T$. A non-conservative, atmospheric-related formulation can be cast in terms of the potential temperature $\theta$ and the Exner function $\pi$ which are defined as
\begin{equation}
\theta\equiv T\left(\frac{p}{p_0}\right)^{-R/c_p},\qquad \pi\equiv \left(\frac{p}{p_0}\right)^{R/c_p},
\end{equation}
where $c_p$ is the specific heat of dry air at constant pressure, $R$ is the gas constant for dry air and $p_0$ is atmospheric pressure at ground level. We can linearize this model upon a decomposition into reference states and perturbations of the form $\Psi(x,z,t)=\overline\Psi(z)+\Psi'(x,z,t)$, and assuming a reference isothermal state ( $\overline T$ constant), we can obtain a rescaled model for the perturbations
 \begin{equation}
\pt  q'+\A\px  q'+\B\pz q'+\S q'=0,
\end{equation}
where
\begin{equation}
\A=\left(\begin{array}{cccc}
U&0& 0& c_s\\
0& U & 0 & 0\\
0 & 0 & U& 0\\
c_s & 0 & 0 & U
\end{array}\right),\;
\S=\left(\begin{array}{cccc}
0&0&0& 0\\
0 & 0 & -\mathcal{N} & -\mathcal{L}\\
0 & \mathcal{N} & 0& 0\\
0 & \mathcal{L} & 0& 0
\end{array}\right),
\end{equation}
and
\begin{equation}
\B=\left(\begin{array}{cccc}
0&0& 0& 0\\
0 & 0 & 0 & c_s\\
0 & 0 & 0& 0\\
0 & c_s & 0 & 0
\end{array}\right),\;
\mathcal{L}=c_s\left[\frac{1}{2\overline\rho}\pz\overline\rho+\frac{1}{\overline\theta}\pz\overline\theta\right].
\end{equation}
We refer the reader to \cite{durran} for the complete derivation of this set of equations.
The unknowns of the presented model are the velocity field perturbation in the x-z plane $(u',w')$ , the potential temperature perturbation $\theta'$ and the Exner function perturbation $\pi'$, all of them represented in the unknown vector $q'=(u',w',\theta',\pi')^t$. Parameters in the isothermal atmosphere are usually initialized as follows. We prescribe a constant value for the reference temperature $\overline{T}$ and the Brunt-V\"ais\"al\"a frequency $\mathcal{N}$. The reference state for the potential temperature is then  given by the relation
\begin{equation}
\mathcal{N}=\sqrt{g\frac{\pz \overline\theta}{\overline\theta}}\,,\qquad \overline\theta(0)=\overline T.
\end{equation}
The reference state for the Exner function is obtained via the hydrostatic balance relation
\begin{equation}
c_p\,\overline\theta\pz\overline\pi=-g\,,\qquad \overline\pi(0)=1.
\end{equation}
The sound speed $c_s$ is constant, and can be computed from
\begin{equation}
c_s=\sqrt{\frac{c_p}{c_v}R\,\overline\pi\,\overline\theta},
\end{equation}
with $c_v$ specific heat of dry air at constant volume, and the reference state for the density is obtained from the state equation of ideal gas
\begin{equation}
\overline\rho=\frac{p_0\overline\pi^{c_v/R}}{R\overline\theta}.
\end{equation}
The remaining parameter $U$ corresponds to the reference state for the horizontal velocity. For all the tests in this section, we set  $\overline{T}=300[K]$, $\mathcal{N}=0.01[s^{-1}]$, $g=9.8[m\,s^{-2}]$, $c_v=717[J\,kg^{-1}\,K^{-1}] $, $c_p=1006 [J\,kg^{-1}\,K^{-1}] $, $R=287.04 [J\,kg^{-1}\,K^{-1}]$, $p_0=10^5[Pa]$ and $U=20[m\,s^{-1}]$.

So far we have presented a method for the resolution of scalar equations; with this example we extend it to constant-coefficient linear hyperbolic systems. As the dimensionally reduced model only treats it as an unidimensional system of balance laws in the x-direction, it is convenient at this point to transform the system into characteristic variables in that direction. Expressing the equations in terms of the characteristic vector $ q$ leads to
\begin{equation}
\pt q+\Lambda \px q+\tilde \B\pz q+\tilde \S q=0,
\end{equation}

where

\begin{align}
\Lambda = Q ^{-1} \A Q = \left(\begin{array}{cccc} U+c_s & 0 &0 &0\\
0 & U-c_s & 0 &0\\
0 & 0 & U & 0\\
0&0&0&U \end{array}\right),\;\tilde\B =Q^{-1}\B Q =\left(\begin{array}{cccc} 0 & 0 &0 &\frac12 c_s\\
0 & 0 & 0 &\frac12 c_s\\
0 & 0 & 0 & 0\\
c_s&c_s&0&0\end{array}\right),
\end{align}
and
\begin{align}
\tilde\S = Q^{-1}\S Q=\left(\begin{array}{cccc} 0 & 0 & 0 &\frac12 \mathcal{L}\\
0 & 0 & 0 &\frac12\mathcal{L}\\
0 & 0 & 0& \mathcal{N}\\
-\mathcal{L}& -\mathcal{L}&-\mathcal{N} & 0\end{array}\right),\; q =Q^{-1} q' =\left[\begin{array}{c} \frac12(\pi'+u')\\ \frac12(\pi'-u')\\ \theta' \\ w'\end{array}\right].
\end{align}
For this example we consider three different vertical discretizations: piecewise linear continuous finite elements, and piecewise constant and linear discontinuous Galerkin expansions. We close the system with periodic boundary conditions at $x=L_x$ for every variable, and we set $w'=0$ at $z=0$ and $z=H_z$.

In the case of the piecewise linear continuous Galerkin method, the extension to the system is straightforward; after uniformly partitioning the vertical interval into $N$ elements, every variable is expanded as a linear combination of the same basis functions, and testing every equation separately leads to an assembled representation of the form:
\begin{equation}
\pt \vq_h + \mathbf{\Lambda}\px \vq_h +\mathbf{C}\, \vq_h +\mathbf{S}\, \vq_h=0,
\end{equation}
where
\begin{align}
\mathbf{\Lambda} = \left(\begin{array}{cccc}
(U+c_s)\mathbf{I} & \mathbf{0} & \mathbf{0}  &\mathbf{0}\\
\mathbf{0} & (U-c_s)\mathbf{I}  & \mathbf{0}  &\mathbf{0}\\
\mathbf{0}        & \mathbf{0} & U\mathbf{I} & \mathbf{0}\\
\mathbf{0}        & \mathbf{0} &  \mathbf{0} & U\mathbf{I}\end{array}\right),\;
\mathbf{S}=\frac12\left(\begin{array}{cccc} \mathbf{0} & \mathbf{0} &\mathbf{0} &\mathcal{L}\mathbf{I}\\
\mathbf{0} & \mathbf{0} & \mathbf{0} &\mathcal{L}\mathbf{I}\\
\mathbf{0} &\mathbf{0} &\mathbf{0}& 2\mathcal{N}\mathbf{I}\\
-2\mathcal{L}\mathbf{I}& -2\mathcal{L}\mathbf{I}&-2\mathcal{N}\mathbf{I} & 0\end{array}\right),
\end{align}
and
\begin{align}
\mathbf{C} =\frac{c_s}{2}\left(\begin{array}{cccc} \mathbf{0} & \mathbf{0} &\mathbf{0} &\mathbf{M}^{-1}\mathbf{K}\\
\mathbf{0} & \mathbf{0} & \mathbf{0} &\mathbf{M}^{-1}\mathbf{K}\\
\mathbf{0} & \mathbf{0} & \mathbf{0} & \mathbf{0}\\
2\mathbf{M}^{-1}\mathbf{K}&2\mathbf{M}^{-1}\mathbf{K}&\mathbf{0}&\mathbf{0}\end{array}\right),\;
\vq_h= \left[\begin{array}{c}q^1_0\\ \vdots \\ q^1_{N} \\ \vdots \\ q^4_0 \\ \vdots \\q^4_{N}\end{array}\right],
\end{align}
with $\mathbf{I}$ the identity matrix of size $N\times N$, and $\mathbf{M},\,\mathbf{K}$ are defined as in eqns. (\ref{matricesfem}).

For piecewise linear discontinuous Galerkin method (we omit the piecewise constant case as it is trivially contained) the only additional consideration is related with the definition of the upwind flux for a linear system. For the $i-th$ element the local equations read:
\begin{align}
\pt c_{0,i}&+\Lambda\px c_{0,i}+\tilde\S\, c_{0,i}+\frac{1}{\dz}\,\{\,\tilde\B\, q(z=z_i)-\tilde\B\, q(z=z_{i-1})\,\}=0,\\
\pt c_{1,i}&+\Lambda\px c_{1,i}+\tilde\S\, c_{1,i}+\frac{3}{\dz}\,\{\,\tilde\B\, q(z=z_i)+\tilde\B\, q(z=z_{i-1})-2\tilde\B\, c_{0,i}\,\}=0.
\end{align}
In this case, being $\tilde\B$ a constant-coefficient matrix that admits a diagonal representation $\Lambda_b$ such that $\tilde\B=R^{-1}\Lambda_b R$, we split $\Lambda_b$ into $\Lambda_b^+$ and $\Lambda_b^-$, containing only non-negative and non-positive eigenvalues respectively. We define $\tilde\B^+=R^{-1}\Lambda_b^+R$ and $\tilde\B^-=R^{-1}\Lambda_b^-R$, and therefore
the upwinding is defined by
\begin{align}
\tilde\B\, q(z=z_{i-1})&= \tilde\B^+\,q_{i-1}(\xi=1)+\tilde\B^-\,q_i(\xi=-1),\\
\tilde\B\,q(z=z_{i-1})&= \tilde\B^+\,q_i(\xi=1)+\tilde\B^-\,q_{i+1}(\xi=-1),
\end{align}
with
\begin{align}
\tilde\B^-=\frac{c_s}{4}\left(\begin{array}{cccc} -1 & -1 &0 &1\\
-1 & -1 &0 &1\\
0 & 0 &0 &0\\
2 & 2&0 &-2\end{array}\right),\quad
\tilde\B^+=\frac{c_s}{4}\left(\begin{array}{cccc}1 & 1 &0 &1\\
1 & 1 &0 &1\\
0 & 0 &0& 0\\
2 & 2 &0 &2\end{array}\right).
\end{align}

\subsubsection{Convergence rates}
 Our first test is performed in order to validate the theoretical convergence rates of the ADER$m$ scheme for hyperbolic balance laws in the presence of a stable source term. For this purpose we set an initial  potential temperature perturbation \begin{equation}\theta'(x,z,0)=\sin\left(\frac{10\pi}{L_x}x\right)\sin\left(\frac{2\pi}{H_z}z\right)\end{equation} and every other perturbation is initialized at rest; values of other test parameters are $L_x=300[km]$ and $H_z=10[km]$. We use a CFL number $CFL=0.5$, and the time step is limited by the fastest wave of the model, i.e., \begin{equation}\dt=CFL\left(\frac{\dx}{U+c_s}\right),\end{equation} being all the simulations stable in this case. Observation time is at $t=500 [s]$. We consider a fixed number of 20 vertical elements. Results for relative errors and convergence rates in space at final time of simulation, with respect to a reference solution with 1024 horizontal elements, for FEM, DGP0 and DGP1 discretizations are shown in Tables \ref{t3}-\ref{t4}. We note that, excepting for the scheme with source term arising from a DGP0 discretization, convergence rates in $L_1$ and $L_{\infty}$ norm, for spatial order 2, 3 and 4 are in accordance to the expected values, i.e., the presence of the source terms arising from vertical discretizations is not a drawback to recover the same orders as in the purely advective problem. Henceforth further numerical studies are performed only with optimally convergent schemes, i.e., with FEM and DGP1 source terms.

\begin{table}[H]
\centering
\caption{Convergence rates for the ADER scheme with a FEM source term for the atmospheric model.}\label{t3}\vskip 5mm
\begin{tabular}{l c c c c c}
\hline\\
Method & N & $L_{\infty}$ error &  $L_{\infty}$ order & $L_1$ error &  $L_1$ order \\ [0.5ex]
\hline\\

ADER2-FEM & 32 & 1.2854e-002 &  &1.3350e-002& \\
& 64 &   2.3160e-002& --- & 2.2032e-002&--- \\
& 128 &  4.9002e-003   & 2.2407  &4.5064e-003 &2.2896\\
& 256 &   1.3337e-003  &1.8774  &  1.1366e-003&1.9872  \\
& 512 &    3.7896e-004 &1.8154 &2.9054e-004 & 1.9679\\ \\

ADER3-FEM & 32 & 3.0623e-002 &  &  3.0548e-002& \\
& 64 &    3.5159e-003 & 3.1227&3.5144e-003&3.1197  \\
& 128 &   3.5367e-004 & 3.3134 &  3.5363e-004& 3.3130\\
& 256 & 3.3499e-005 &3.4002  & 3.3502e-005  & 3.3999 \\
& 512 &  3.4205e-006 & 3.2919  &3.4204e-006& 3.2920\\ \\

ADER4-FEM & 32 & 2.9463e-003 &  &2.9301e-003& \\
& 64 &  4.0365e-004 & 2.8677 &4.0416e-004&2.8580 \\
& 128 & 4.3203e-005 &3.2239 &4.3196e-005 &3.2259\\
& 256 & 3.0715e-006 & 3.8141  &3.0719e-006&3.8137\\
& 512 &1.8900e-007 & 4.0225 & 1.8899e-007&4.0227\\

\end{tabular}
\end{table}

\begin{table}
\centering
\caption{Convergence rates for the ADER scheme with a DGP0 and DGP1 source terms for the atmospheric model.}\label{t4}\vskip 8mm
\begin{tabular}{l c c c c c}
\hline\\
Method & N & $L_{\infty}$ error &  $L_{\infty}$ order & $L_1$ error &  $L_1$ order \\ [0.5ex]
\hline\\\\

ADER2-DGP0 & 32 &  1.0842e-001&  &6.3604e-002& \\
& 64 &   4.9767e-002& 1.1233 & 2.6742e-002&1.2500 \\
& 128 &  2.1215e-002& 1.2301 & 8.2006e-003&1.7053\\
& 256 &  9.1311e-003& 1.2163  &2.1918e-003&1.9036 \\
& 512 &  3.7758e-003& 1.2740 &5.6313e-004&1.9606\\ \\

ADER3-DGP0& 32 & 3.4761e-002 &  &2.4806e-002& \\
& 64 &     5.4393e-003&2.6760 &5.4400e-003&2.1890  \\
& 128 &  1.7873e-003&1.6056&  1.7865e-003&1.6065\\
& 256 &6.9451e-004&  1.3637& 6.9449e-004&1.3631 \\
& 512 &2.2774e-004& 1.6086  &2.2774e-004&1.6086\\ \\

ADER4-DGP0& 32 & 1.9014e-002&  &1.1660e-002& \\
& 64 & 3.4347e-003& 2.4688 & 3.0461e-003&1.9366 \\
& 128 &1.4788e-003&1.2158  &1.4671e-003&1.0540\\
& 256 &6.5413e-004&1.1768  &6.5380e-004&1.1660\\
& 512 &2.2264e-004&1.5549 &2.2263e-004& 1.5542\\\\\\

ADER2-DGP1 & 32 &  1.3685e-002 &  &6.6063e-003& \\
& 64 &   6.5137e-003   &1.0710 & 1.9915e-003&1.7300 \\
& 128 &  2.8117e-003   &1.2120  &5.5509e-004 &1.8431\\
& 256 &   1.2360e-003  &1.1858  &  1.4800e-004& 1.9071  \\
& 512 &    5.2020e-004 &1.2485 &3.9667e-005 & 1.8996\\ \\

ADER3-DGP1 & 32 & 8.2739e-004 &  &  8.4437e-004& \\
& 64 &    1.0738e-004 & 2.9458&1.0832e-004&2.9626  \\
& 128 &   1.3539e-005 & 2.9876 &  1.3607e-005& 2.9929\\
& 256 & 1.6960e-006 &2.9969  & 1.7033e-006 &2.9979 \\
& 512 &  2.1211e-007 & 2.9993  &2.1298e-007& 2.9995\\ \\

ADER4-DGP1 & 32 & 2.7928e-004 &  &1.4522e-004& \\
& 64 &  8.7557e-006 & 4.9953 &4.8145e-006&4.9147 \\
& 128 &  2.6248e-007 &5.0599 & 2.2554e-007& 4.4160\\
& 256 & 1.2429e-008 &4.4004  &1.2793e-008& 4.1400\\
& 512 &7.2632e-010 & 4.0970 &7.3387e-010&4.1236\\

\end{tabular}
\end{table}

\subsubsection{Inertia-gravity waves in a periodic channel}

Having verified the numerical accuracy of the scheme, we turn our attention to a test case which is meaningful from an atmospheric modelling point of view. We study the behavior of non-hydrostatic scale gravity waves (NHIGW). In the same framework of \cite{skamarock} and \cite{testjcp}, with the same parameters as in the previous section, we test with an initial potential temperature perturbation given by:
\begin{equation}
\theta'(x,z,0)=0.01\frac{\sin\left(\frac{\pi z}{H_z}\right)}{1+\left(\frac{x-L_x/3}{a}\right)^2}
\end{equation}

Computational domain and other parameters are given in Table \ref{param}. This test is performed in order to verify proper speed propagation of the inertia-gravity waves and its symmetry with respect to the advected initial perturbation. The initial perturbation is expected to be translated by the background horizontal velocity, while at the same time breaking and radiating symmetrically to the left and to the right; the axis of symmetry should recover the background translation speed.

\begin{table}[H]
\centering
\caption{Computational domain and parameters for simulations of NHIGW test case}\label{param}\vskip 5mm
\begin{tabular}{l c c c c c c c  c}
\hline\\
Test &  $L_x$ & $H_z$ & $a$& $T_f$ & CFL & $\dx$ & $\dz$& $\dt$ \\ [0.5ex]
\hline\\
NHIGW &  $300 [km]$ & $10 [km]$ &$5 [km]$ & $3\times 10^3 [s]$ & 0.5 & $1 [km]$ & $500 [m]$& $1.36 [s]$ \\
\end{tabular}
\end{table}

\begin{figure}
\caption{1. Initial potential temperature perturbation for the NHIGW test case (top) 2. Evolution of $u'^2+w'^2$ during the NHIGW test case (bottom).}\label{nonhyd1}
\centering
\includegraphics[height=.32\textheight,width=\textwidth]{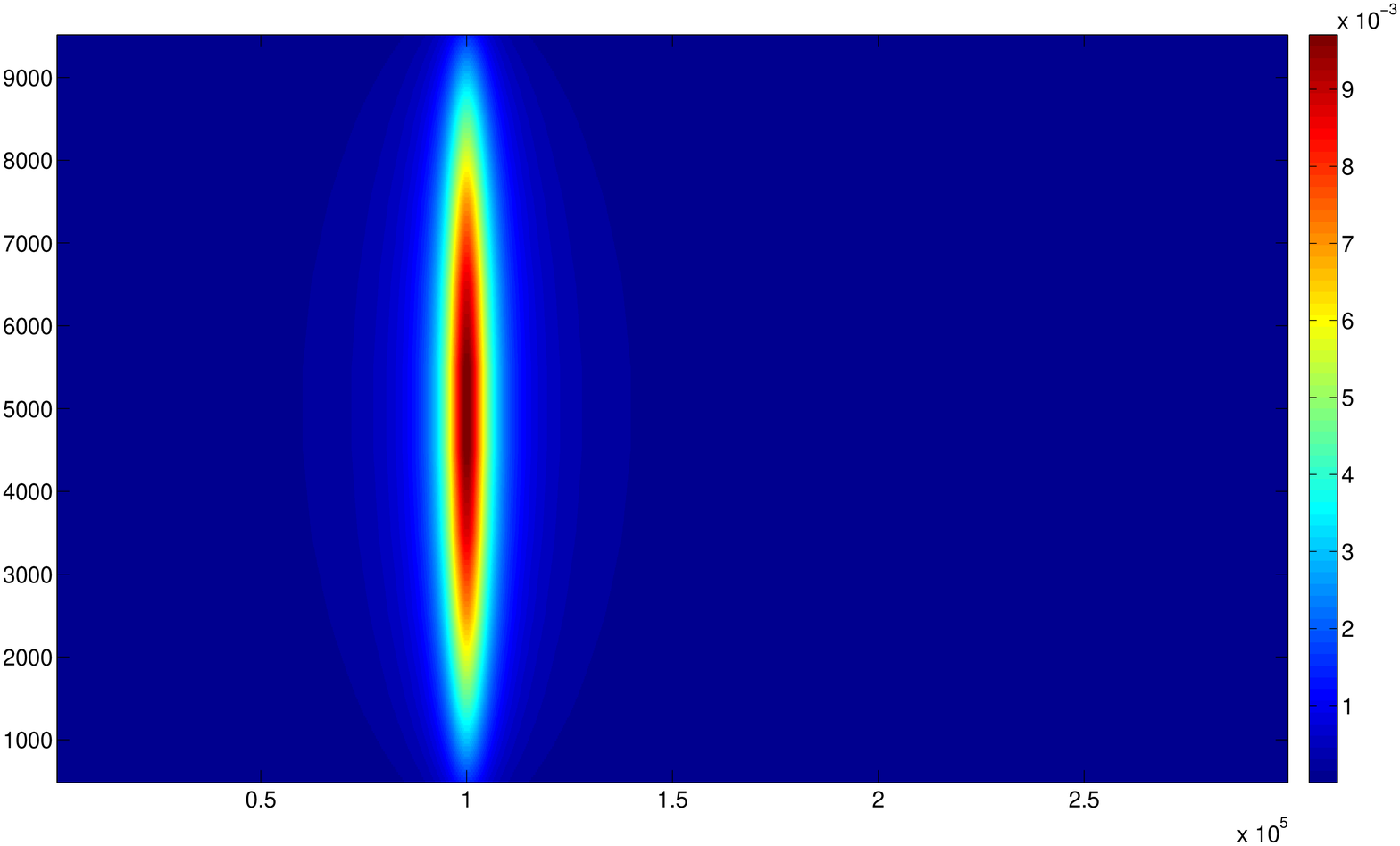}

\includegraphics[height=.32\textheight,width=\textwidth]{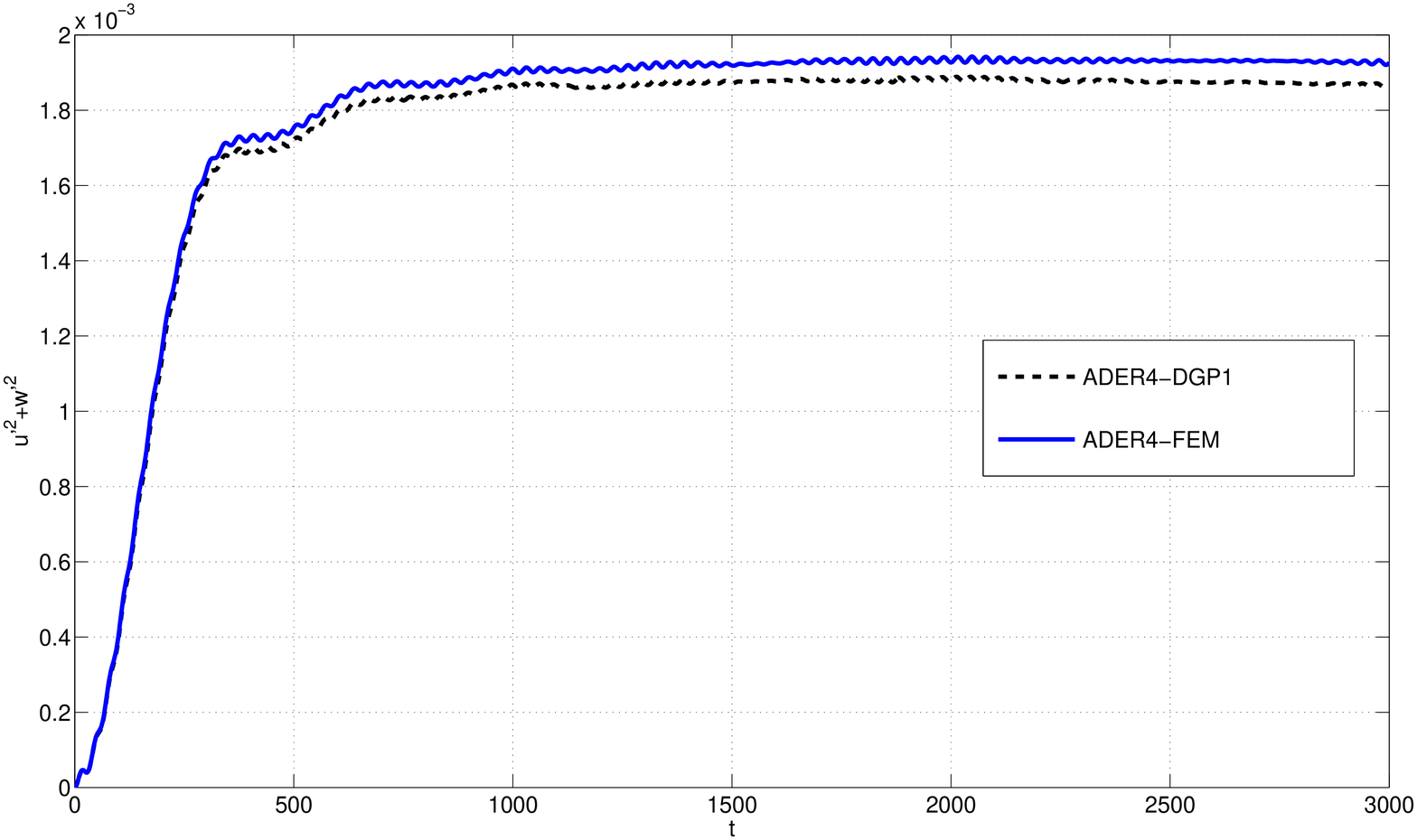}
\end{figure}

\begin{figure}
\caption{ADER4-FEM scheme, NHIGW test, potential temperature perturbation at $z=5,000[m]$  1. At $t=1,000[s]$ (top) 2. At $t=2000[s]$ (middle) 3. At $t=3,000[s]$ (bottom).}\label{nonhyd2}
\centering

\includegraphics[height=.32\textheight,width=\textwidth]{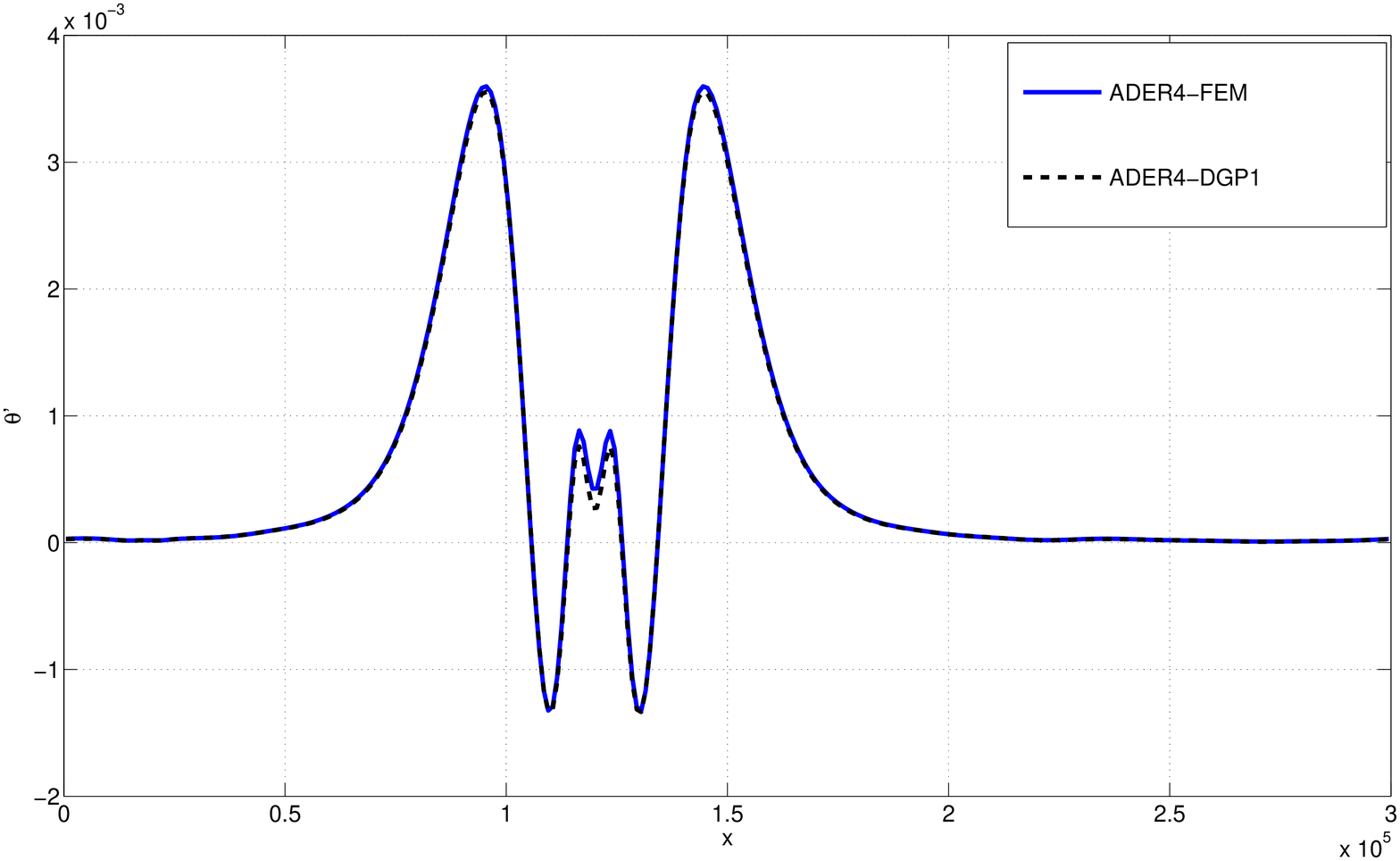}

\includegraphics[height=.32\textheight,width=\textwidth]{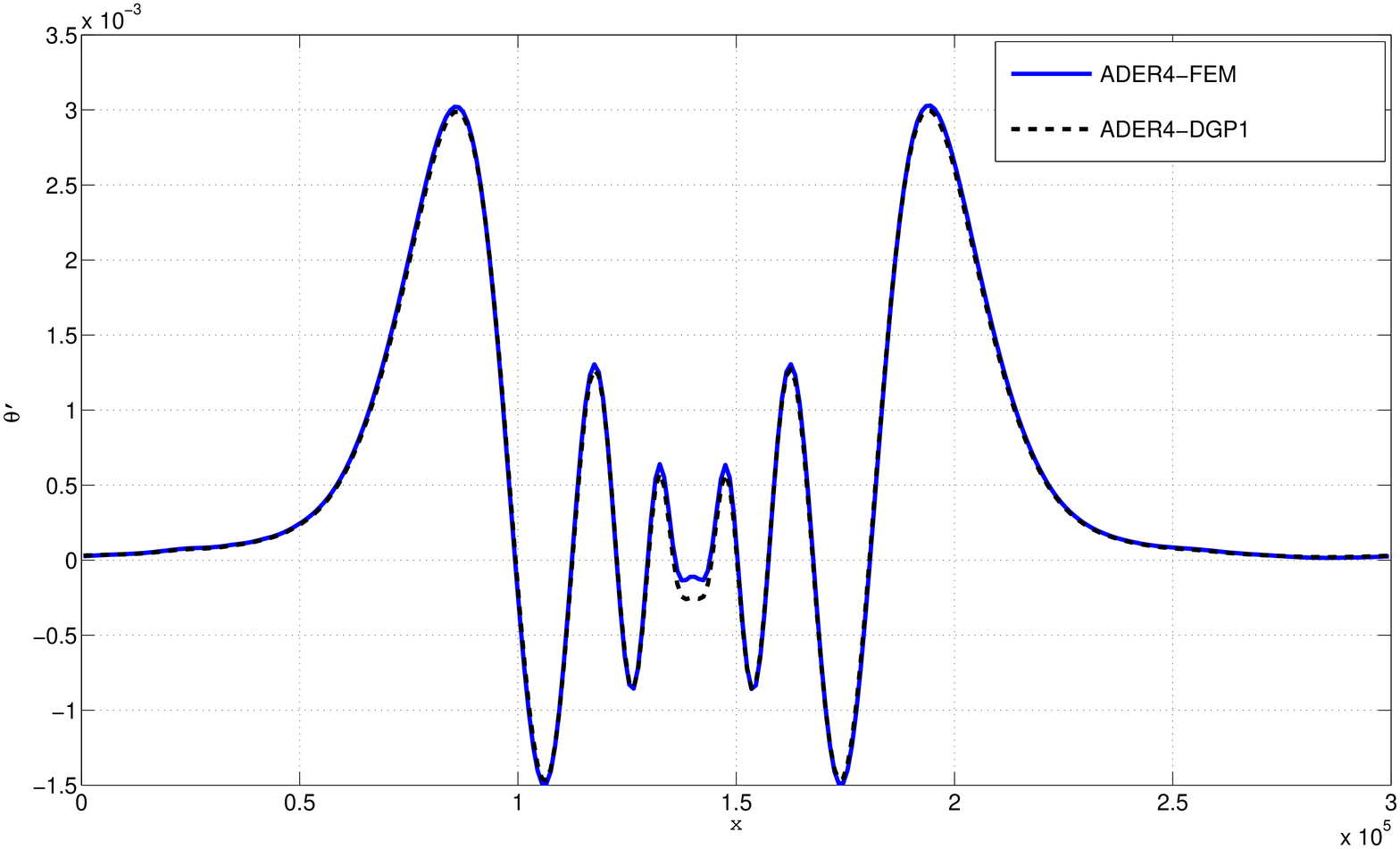}

\includegraphics[height=.32\textheight,width=\textwidth]{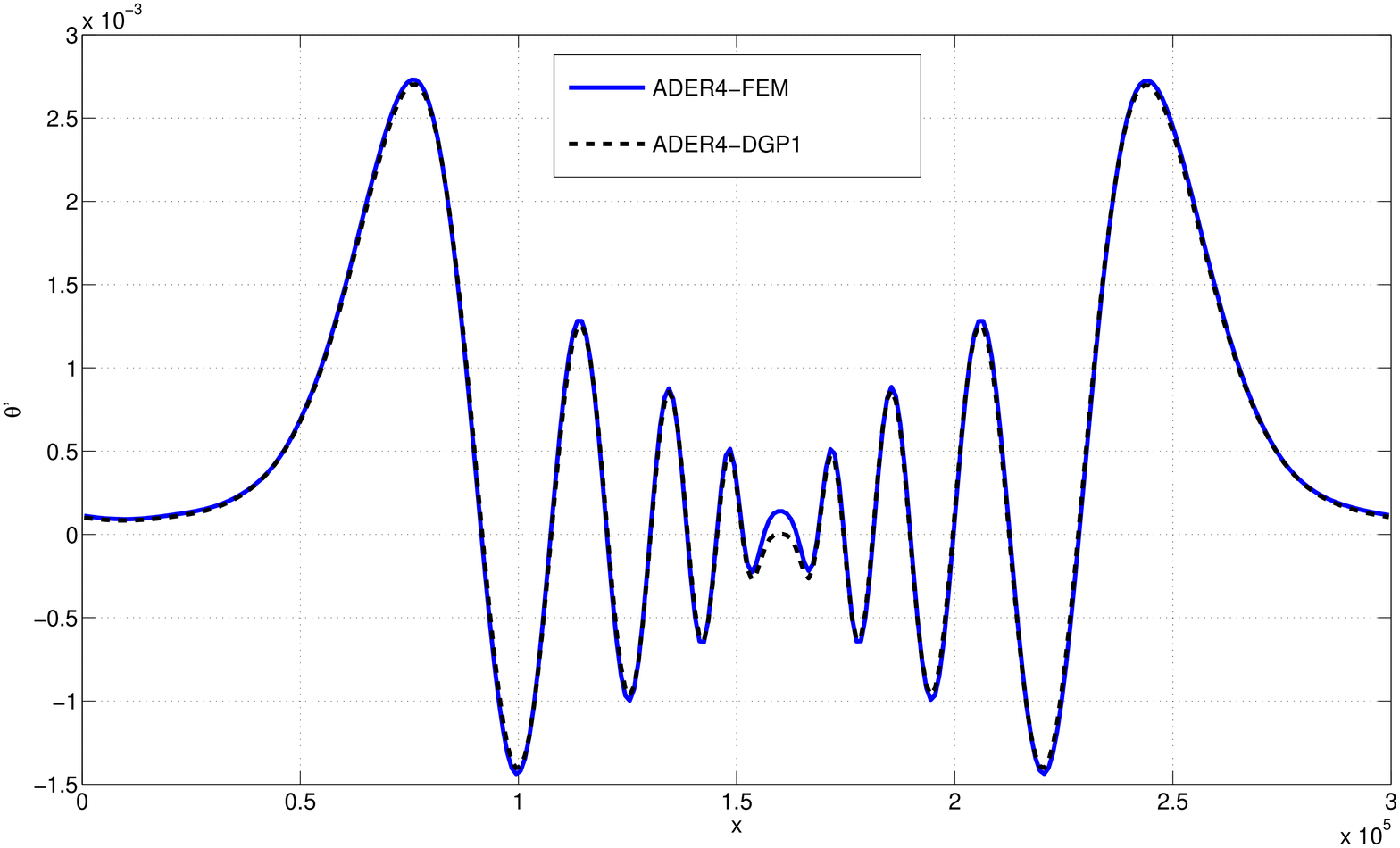}
\end{figure}

\begin{figure}
\caption{ADER4-FEM scheme, NHIGW test, potential temperature perturbation 1. At $t=1,000[s]$ (top) 2. At $t=2,000[s]$ (middle) 3. At $t=3,000[s]$ (bottom).}\label{nonhyd3}
\centering
\includegraphics[height=.32\textheight,width=\textwidth]{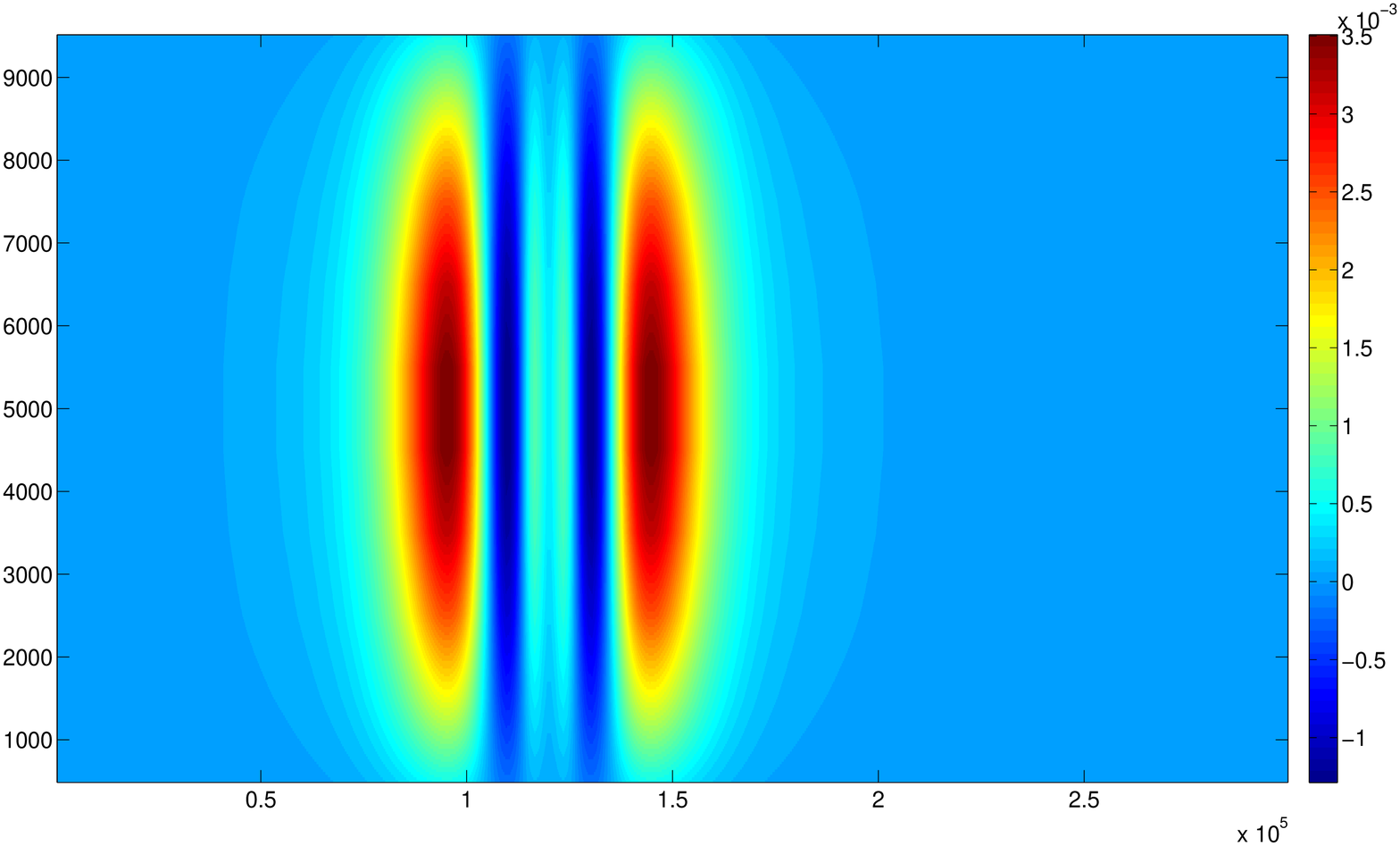}

\includegraphics[height=.32\textheight,width=\textwidth]{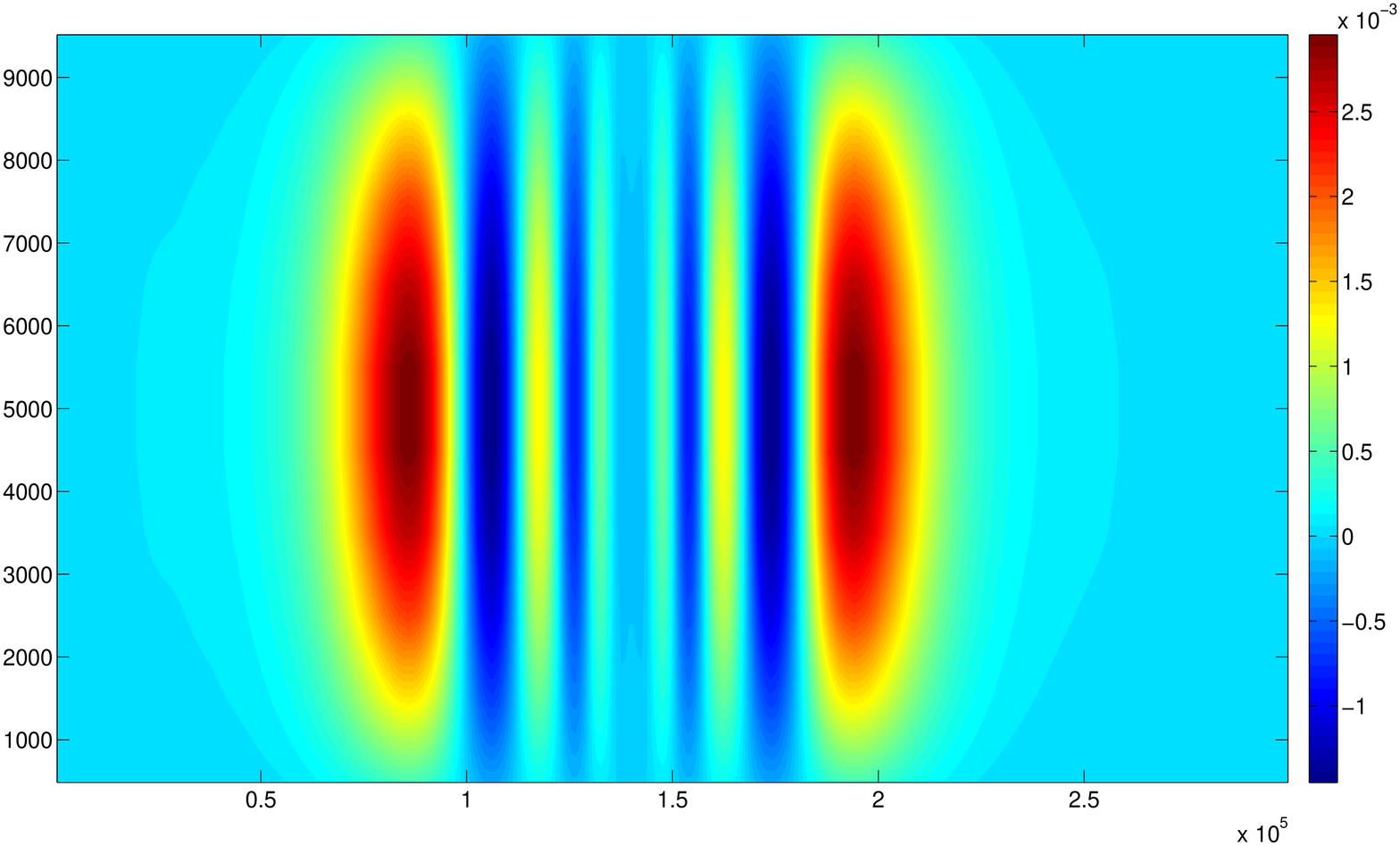}

\includegraphics[height=.32\textheight,width=\textwidth]{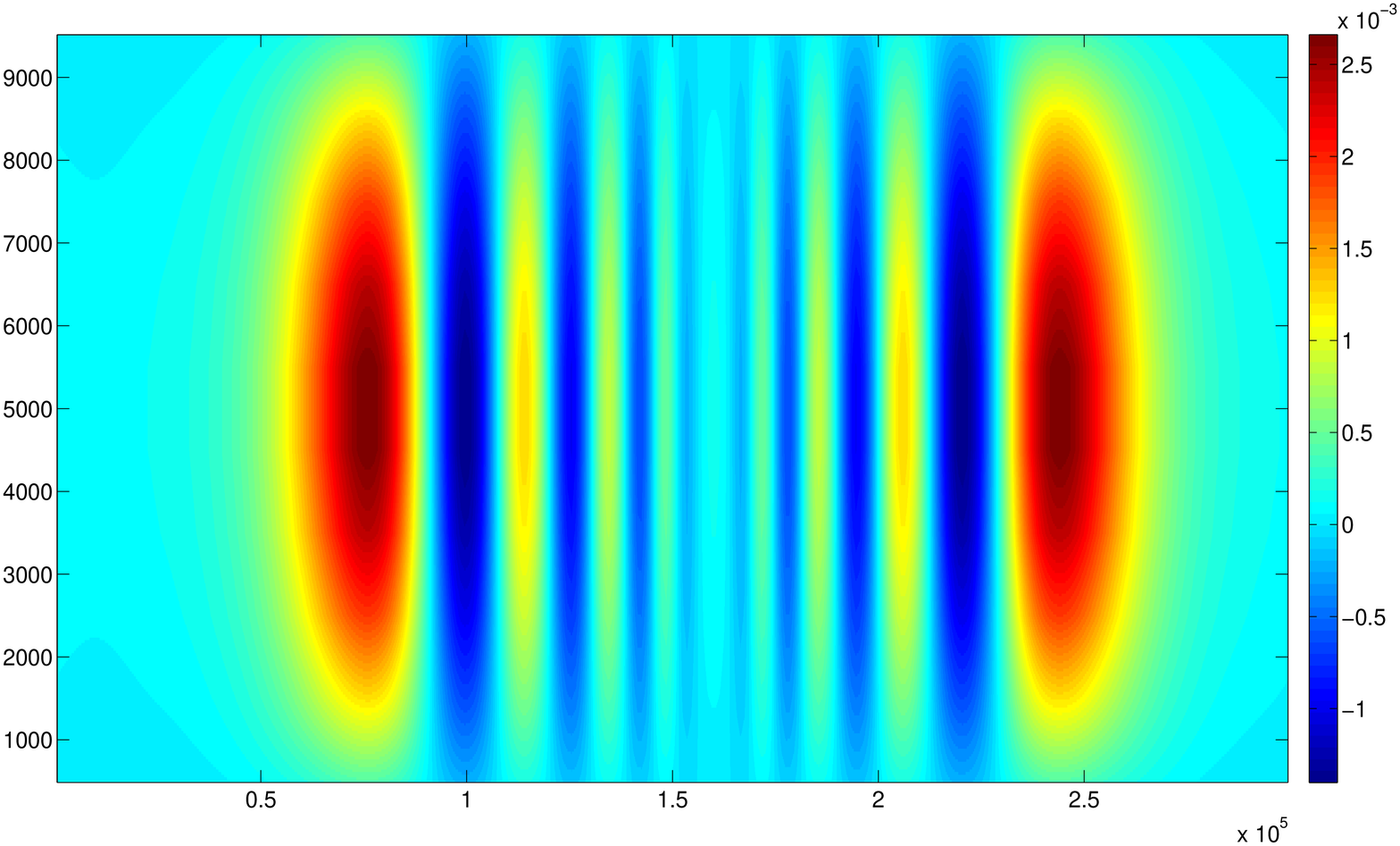}

\end{figure}

\begin{figure}
\caption{ADER4-DGP1 scheme, NHIGW test, velocity field 1. At $t=1,000[s]$ (top) 2. At $t=2,000[s]$ (middle) 3. At $t=3,000[s]$ (bottom).}\label{nonhyd4}
\centering
\includegraphics[height=.32\textheight,width=\textwidth]{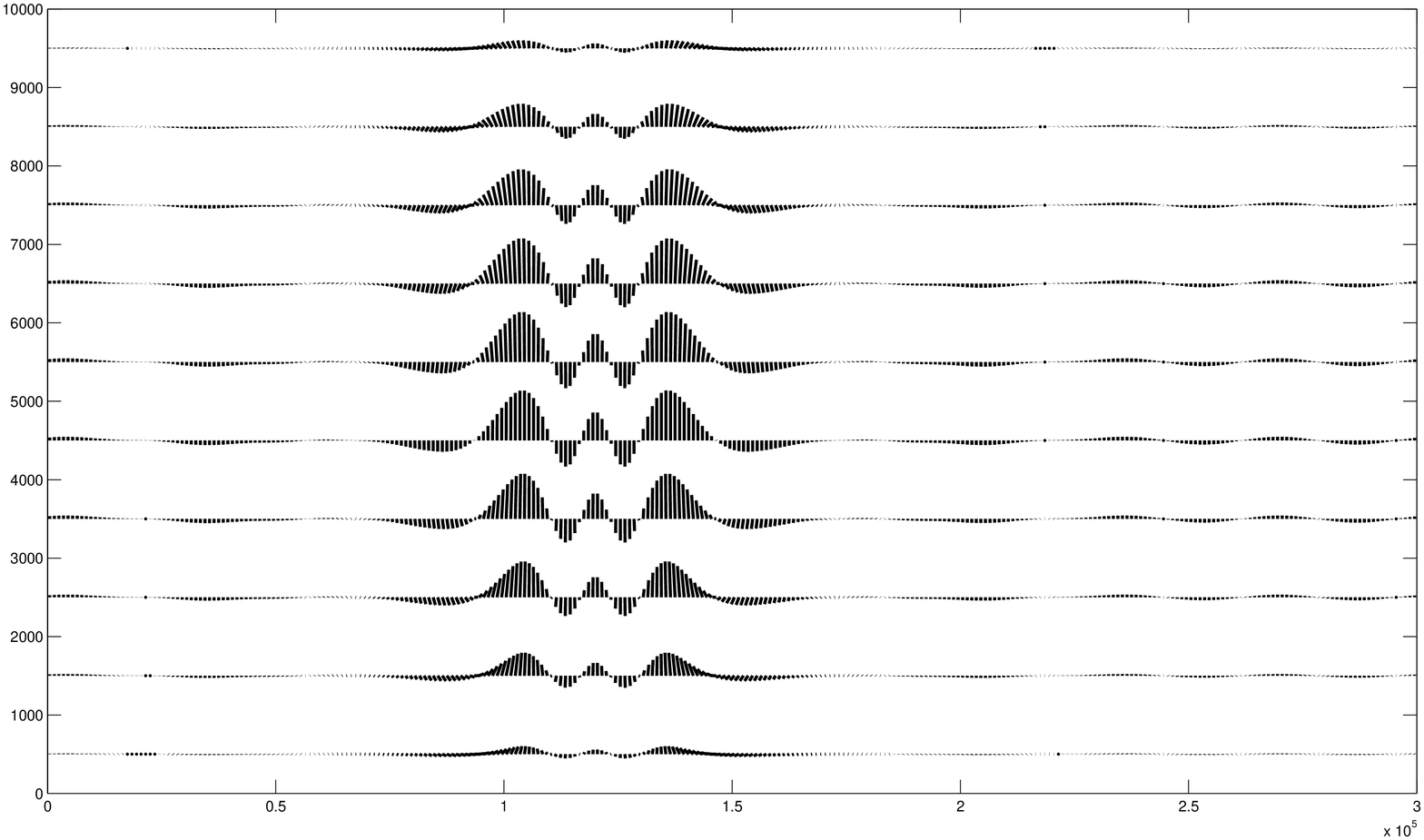}

\includegraphics[height=.32\textheight,width=\textwidth]{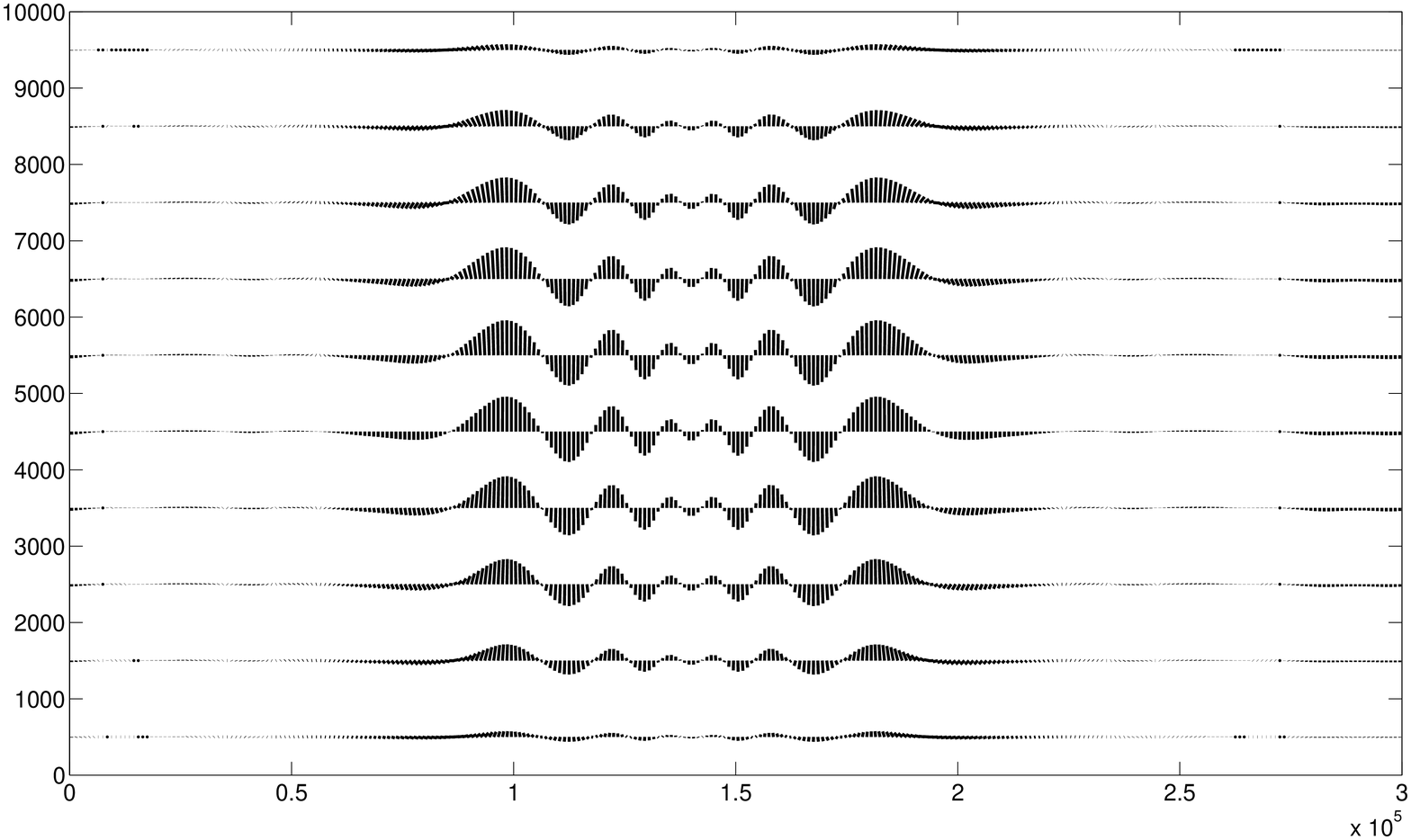}

\includegraphics[height=.32\textheight,width=\textwidth]{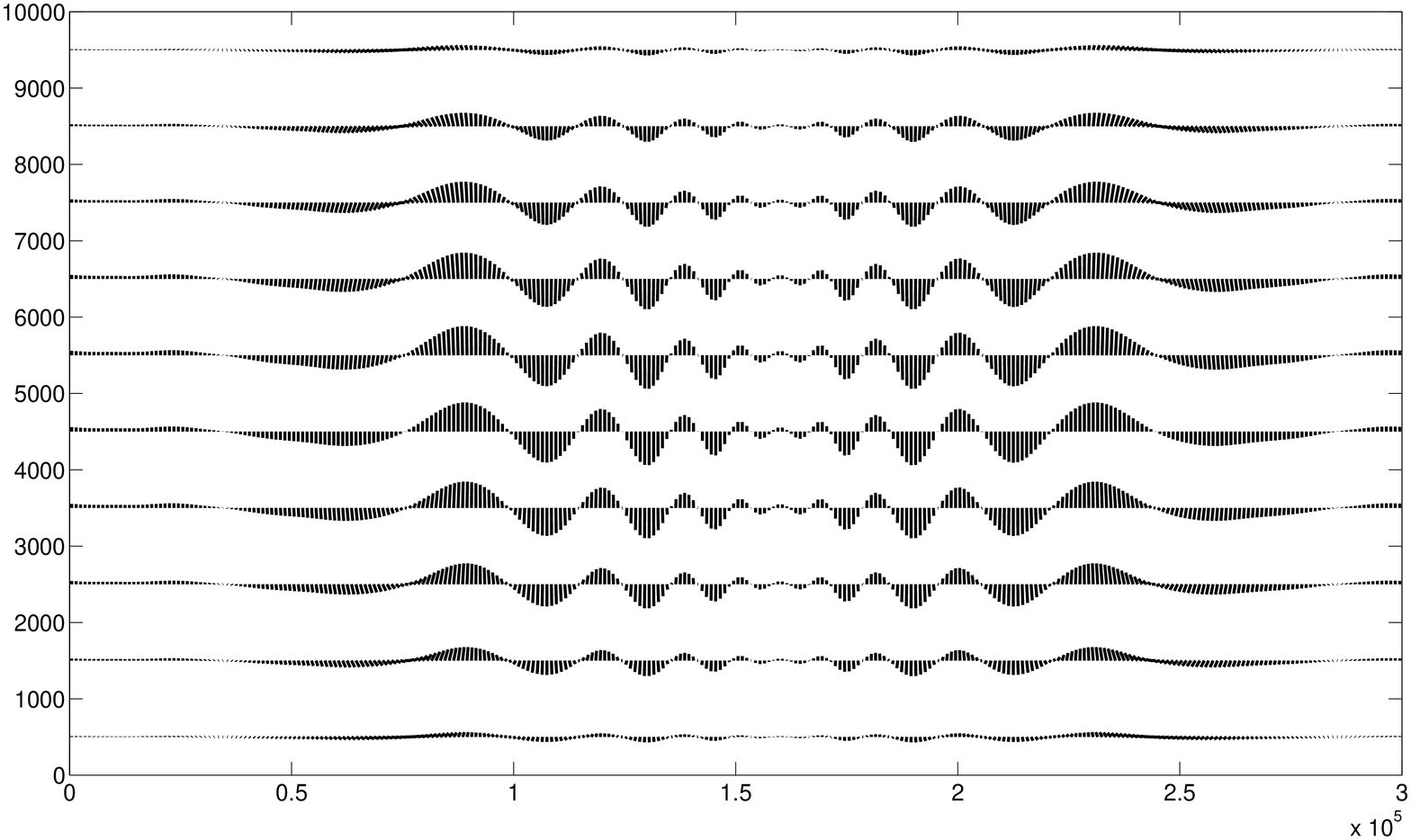}
\end{figure}

\begin{figure}
\caption{ADER4-DGP1 scheme, NHIGW test, divergence of the velocity field 1. At $t=1,000[s]$ (top) 2. At $t=2,000[s]$ (middle) 3. At $t=3,000[s]$ (bottom).}\label{nonhyd5}
\centering
\includegraphics[height=.32\textheight,width=\textwidth]{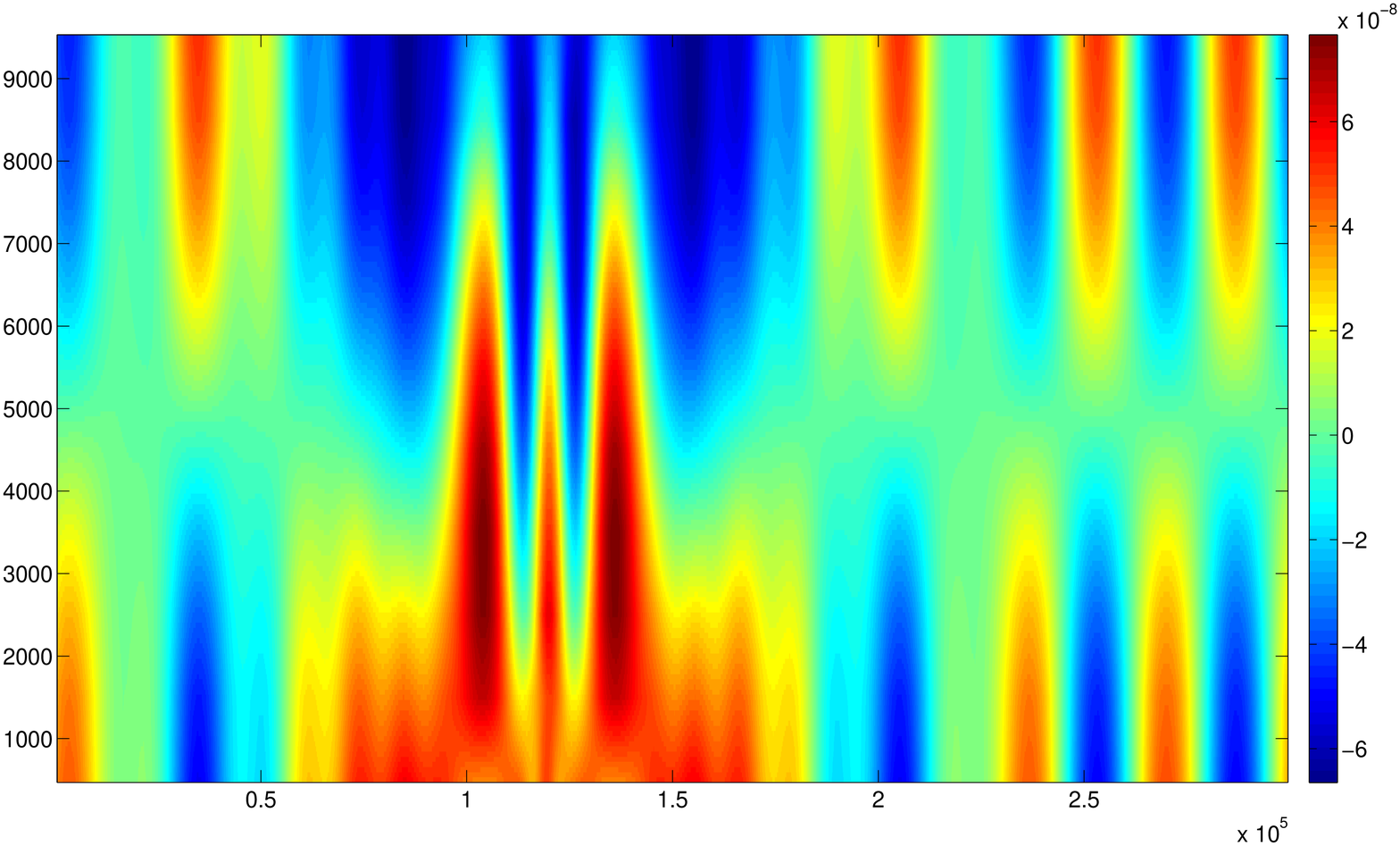}

\includegraphics[height=.32\textheight,width=\textwidth]{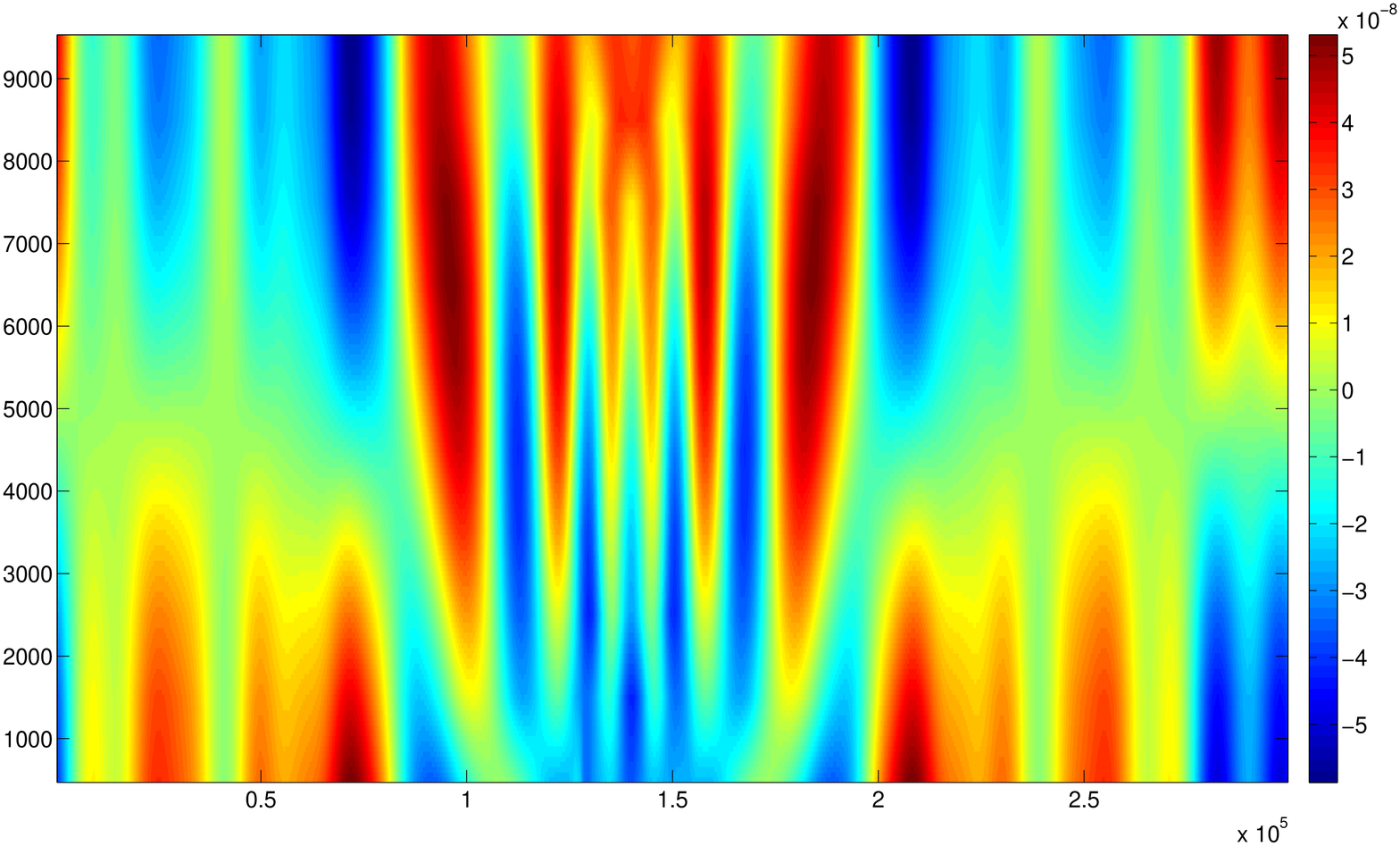}

\includegraphics[height=.32\textheight,width=\textwidth]{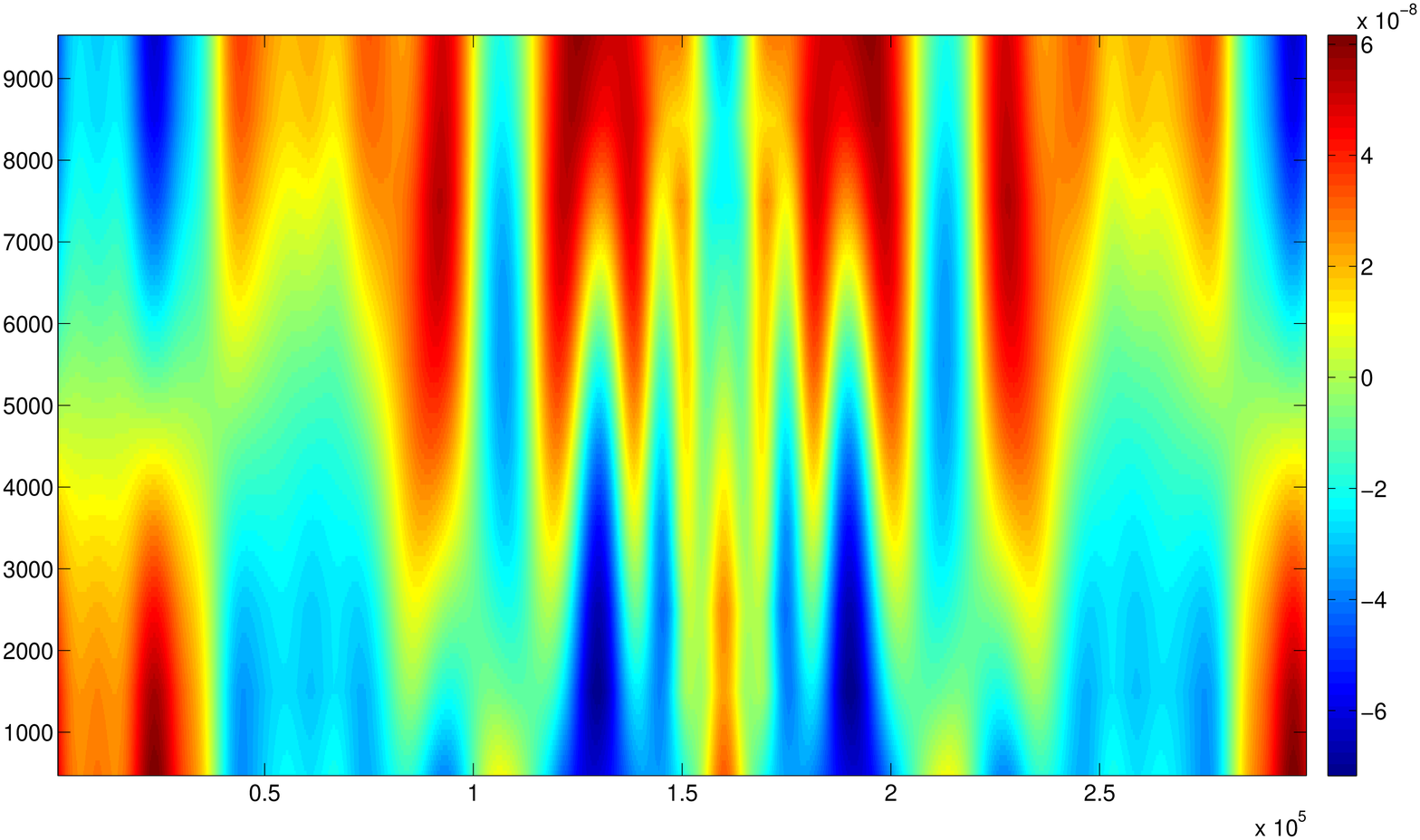}

\end{figure}

As in the convergence rates validation, we first observe that for every method a reduction of the CFL number to .5 is needed in order to have stable simulations. We conjecture that this is due to the presence of a source term which includes an oscillatory behavior that is reflected in a set of purely imaginary eigenvalues on its spectrum, independently of the nature of the vertical discretization. Although nonlinear, the treatment that the ADER scheme gives to the source term can be interpreted as a single-step, explicit method and therefore a reduction of the time step is required in order to include a section the imaginary axis inside the stability region of the scheme.
Table \ref{vv} shows values for the potential temperature extrema that are  in accordance with the results obtained for the NHIGW test case in \cite{testjcp}. Simulations show that the phenomena is accurately described in physical terms as it preserves symmetry in the propagation of the initial perturbation  and the correct advection speed for both FEM and DGP1 source terms (see \ref{nonhyd2},\ref{nonhyd3} and \ref{nonhyd4}). These second-order vertical schemes manage to conserve the quantity $u'^2+w'^2$ as shown in figure \ref{nonhyd1}, which is not only of physical meaning, but also illustrates the stability of the scheme; finally, the velocity field and its divergence (see \ref{nonhyd4} and \ref{nonhyd5}) are qualitatively correct in the sense that they reflect oscillations for the divergence and ``fan-like'' behavior for the velocity field.

\begin{table}[H]
\centering
\caption{Extrema at $t=3,000 [s]$ for the NHIGW test case with ADER4 scheme for horizontal resolution, FEM and DGP1 schemes for vertical discretization}\label{vv}\vskip 5mm
\begin{tabularx}{.99\textwidth}{l c c c c c c c c }
\hline\\
Source & $\theta'_{max}$ & $\theta'_{min}$ &  $u'_{max}$ & $u'_{min}$ & $w'_{max}$ & $w'_{min}$ & $\pi'_{max}$ & $\pi'_{min}$ \\ [0.5ex]
\hline\\
FEM & 2.70e-3 &-1.43e-3&2.66e-3 &-2.66e-3 & 7.69e-4 & -8.93e-4&5.28e-4 &-5.28e-4 \\
DGP1&2.63e-3 &-1.36e-3&2.59e-3 &-2.58e-3 & 7.18e-4& -8.38e-4&4.54e-4 &-4.33e-4 \\
\end{tabularx}
\end{table}

\begin{figure}
\caption{ADER4-FEM scheme, NHIGW test, Exner function perturbation 1. At $t=1,000[s]$ (top) 2. At $t=2000[s]$ (middle) 3. At $t=3,000[s]$ (bottom).}\label{nonhyd6}
\centering
\includegraphics[height=.32\textheight,width=\textwidth]{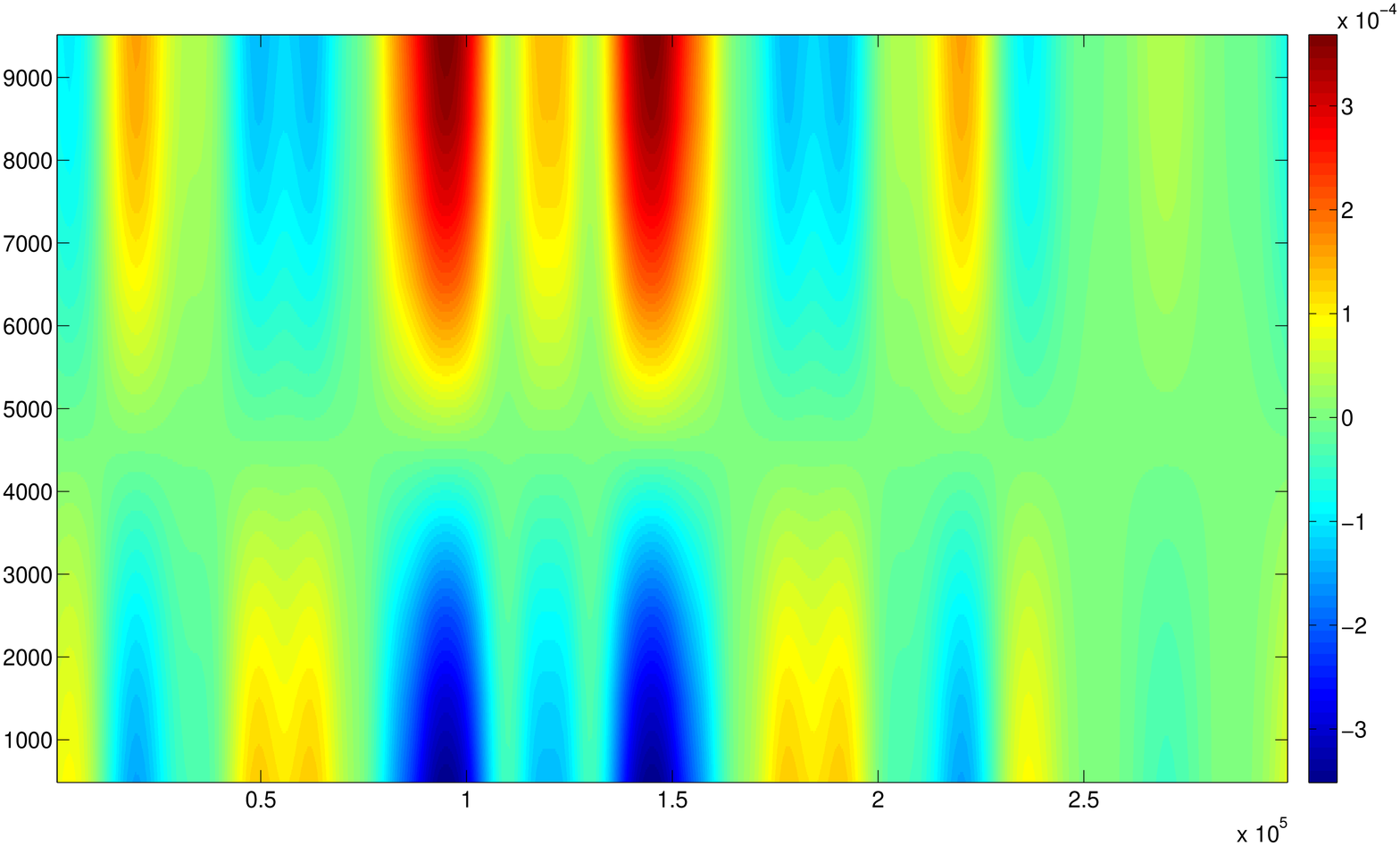}

\includegraphics[height=.32\textheight,width=\textwidth]{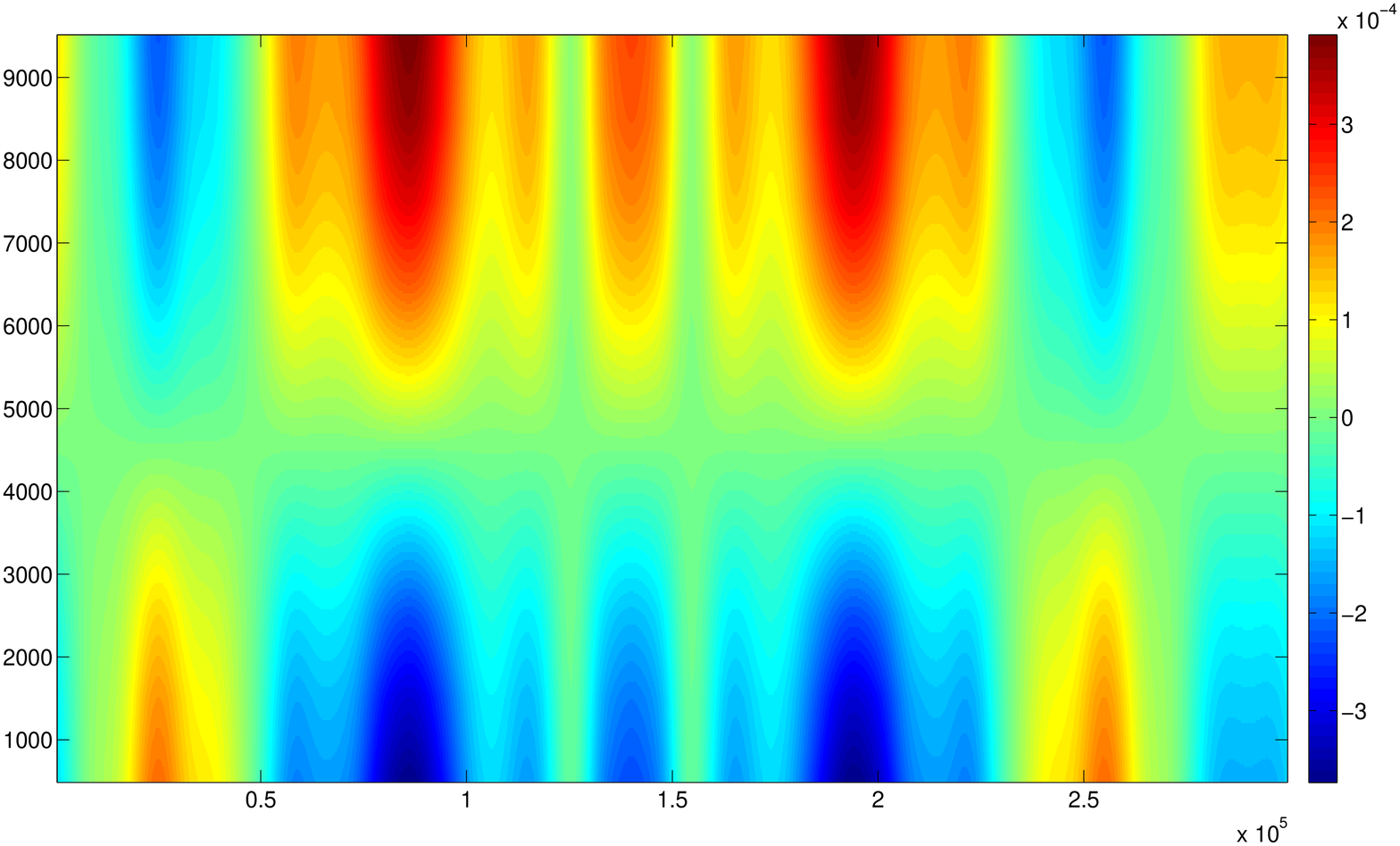}

\includegraphics[height=.32\textheight,width=\textwidth]{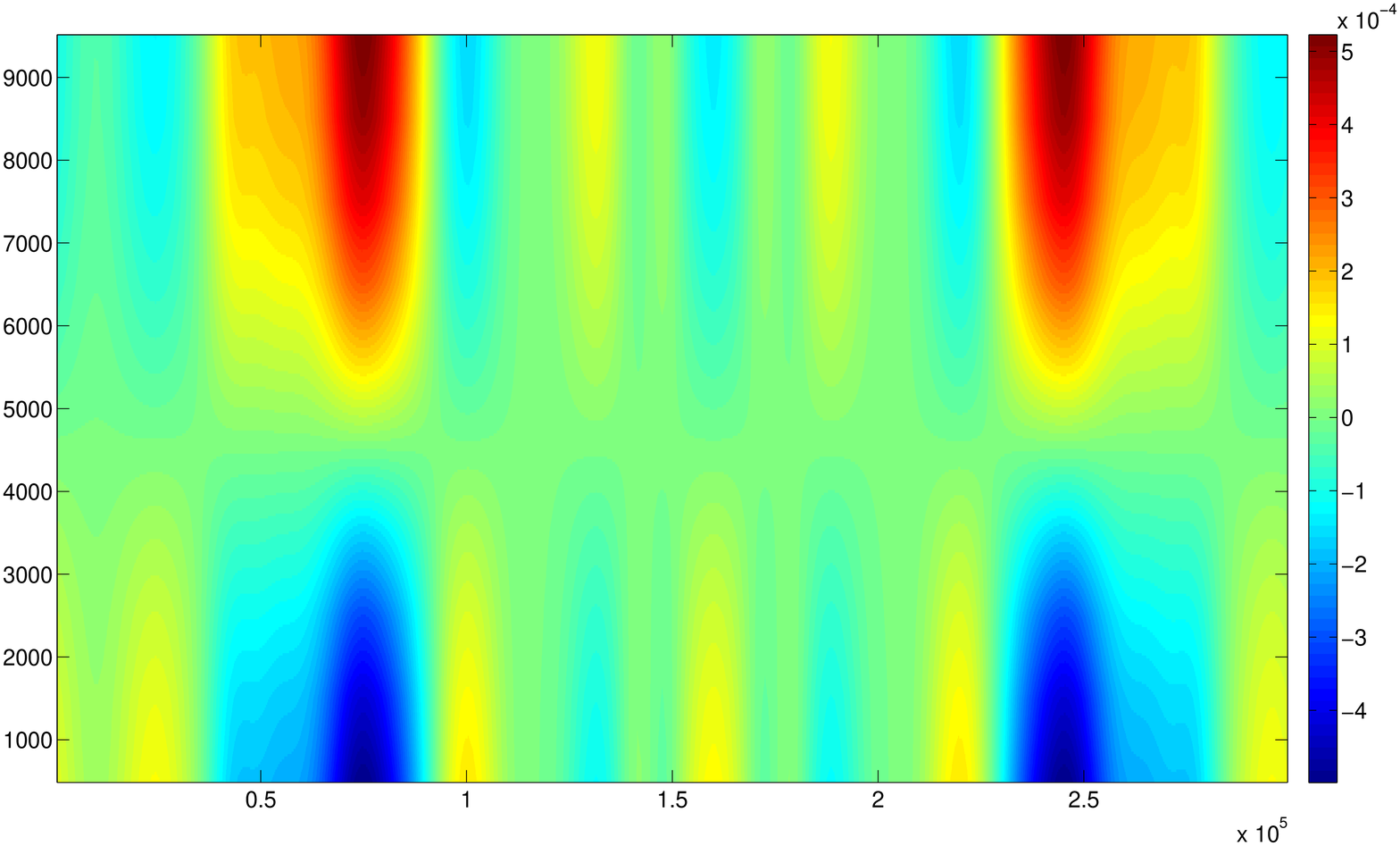}

\end{figure}

\section{Summary and outlook}
We have presented an approach for the generation of layered, high horizontal resolution advective models, which has
been tested in linear advection models. We have implemented different approaches for the layering procedure, while the resulting dimensionally reduced system is solved via the ADER finite volume scheme with prescribed order of accuracy. In the test cases, the method has shown a good performance in terms of stability, convergence, and recovery of the underlying physics of the equations.
A feature of the method is efficiency; the C-K procedure is computed only once before initialization, and the all the structure of the code can be easily parallelized/vectorized. Moreover, the layered structure of the model can be parallelized in the context of domain decomposition. Thus, the extension of the presented procedure to very higher-order approximations seems to be feasible and meaningful.
Although we have limited our analysis to linear model, the discontinuous Galerkin layering approach can be easily extended to non-linear models due to its local character. The class of models where the presented approach can be applied is quite wide, being the main requirement that the advective behavior of the model is governed by the horizontal direction, which is a plausible assumption in most atmospheric phenomena.

\end{document}